\documentclass[]{amsart}
\usepackage[all]{xy}
\usepackage{graphicx}
\usepackage{amscd,amsthm,amssymb,amsfonts,amsmath,euscript}

\theoremstyle{plain}
\newtheorem{thm}{Theorem}[section]
\newtheorem{lm}[thm]{Lemma}
\newtheorem{prop}[thm]{Proposition}
\newtheorem{cor}[thm]{Corollary}

\theoremstyle{definition}

\theoremstyle{remark}
\newtheorem{remark}[thm]{Remark}
\newtheorem*{thank}{Acknowledgments}
\newcommand{\nc}{\newcommand}

\def\makeop#1{\expandafter\def\csname#1\endcsname
  {\mathop{\rm #1}\nolimits}\ignorespaces}
\makeop{Hom}   \makeop{End}   \makeop{Aut}   \makeop{Isom}  \makeop{Pic} 
\makeop{Gal}   \makeop{ord}   \makeop{Char}  \makeop{Div}   \makeop{Lie} 
\makeop{PGL}   \makeop{Corr}  \makeop{PSL}   \makeop{sgn}   \makeop{Spf}
\makeop{Spec}  \makeop{Tr}    \makeop{Nr}    \makeop{Fr}    \makeop{disc}
\makeop{Proj}  \makeop{supp}  \makeop{ker}   \makeop{im}    \makeop{dom}
\makeop{coker} \makeop{Stab}  \makeop{SO}    \makeop{SL}    \makeop{SL}
\makeop{Cl}    \makeop{cond}  \makeop{Br}    \makeop{inv}   \makeop{rank}
\makeop{id}    \makeop{Fil}   \makeop{Frac}  \makeop{GL}    \makeop{SU}
\makeop{Trd}   \makeop{Sp}    \makeop{Tr}    \makeop{Trd}   \makeop{diag}
\makeop{Res}   \makeop{ind}   \makeop{depth} \makeop{Tr}    \makeop{st}
\makeop{Ad}    \makeop{Int}   \makeop{tr}    \makeop{Sym}   \makeop{can}
\makeop{length}\makeop{SO}    \makeop{torsion} \makeop{GSp} \makeop{Ker}
\def\makebb#1{\expandafter\def
  \csname bb#1\endcsname{{\mathbb{#1}}}\ignorespaces}
\def\makebf#1{\expandafter\def\csname bf#1\endcsname{{\bf
      #1}}\ignorespaces} 
\def\makegr#1{\expandafter\def
  \csname gr#1\endcsname{{\mathfrak{#1}}}\ignorespaces}
\def\makescr#1{\expandafter\def
  \csname scr#1\endcsname{{\EuScript{#1}}}\ignorespaces}
\def\makecal#1{\expandafter\def\csname cal#1\endcsname{{\mathcal
      #1}}\ignorespaces} 

\def\doLetters#1{#1A #1B #1C #1D #1E #1F #1G #1H #1I #1J #1K #1L #1M
                 #1N #1O #1P #1Q #1R #1S #1T #1U #1V #1W #1X #1Y #1Z}
\def\doletters#1{#1a #1b #1c #1d #1e #1f #1g #1h #1i #1j #1k #1l #1m
                 #1n #1o #1p #1q #1r #1s #1t #1u #1v #1w #1x #1y #1z}
\doLetters\makebb   \doLetters\makecal  \doLetters\makebf
\doLetters\makescr 
\doletters\makebf   \doLetters\makegr   \doletters\makegr
     \def\qed{\qedmark\medbreak}%
\def\qedmark{{\enspace\vrule height 6pt width 5pt depth 1.5pt}}%

\normalsize

\makeop{Bl}

\def\Fpbar{\overline{\bbF}_p}
\def\Fp{{\bbF}_p}

\def\Qbar{\overline{\bbQ}}

\newcommand{\Z}{\mathbb Z}

\newcommand{\R}{\mathbb R}
\newcommand{\C}{\mathbb C}
\newcommand{\F}{\mathbb F}

\newcommand{\npr}{\noindent }
\newcommand{\pr}{\indent }

\newcommand{\<}{\langle}   
\renewcommand{\>}{\rangle} 

\nc{\embed}{\hookrightarrow}

\newcommand{\ch}{characteristic }
\newcommand{\ac}{algebraically closed }
\newcommand{\dieu}{Dieudonn\'{e} }

\nc{\ol}{\overline}
\nc{\wt}{\widetilde}
\nc{\opp}{\mathrm{opp}}
\def\ul{\underline}

\makeop{Ram}
\makeop{Rep}

\begin{document}
\renewcommand{\thefootnote}{\fnsymbol{footnote}}
\setcounter{footnote}{-1}
\numberwithin{equation}{section}

\title{Irreducibility and $p$-adic monodromies on the Siegel
  moduli spaces}
\author{Chia-Fu Yu}
\address{
Institute of Mathematics \\
Academia Sinica \\
128 Academia Rd.~Sec.~2, Nankang\\ 
Taipei, Taiwan \\ and NCTS (Taipei Office)}
\address{
Max-Planck-Institut f\"ur Mathematik \\
Vivatsgasse 7 \\
Bonn, 53111\\ 
Germany}
\email{chiafu@math.sinica.edu.tw}

\date{February 4, 2008, revised .  The research is partially supported by NSC
 96-2115-M-001-001.}

\begin{abstract}
We generalize the surjectivity result of the $p$-adic monodromy for
the ordinary locus of a Siegel moduli space by Faltings and Chai
(independently by Ekedahl) to that for any $p$-rank stratum. 
We discuss irreducibility and connectedness of some $p$-rank 
strata of the moduli spaces with parahoric level structure. Finer
results are obtained on the Siegel 3-fold with Iwahori level
structure.      
\end{abstract}
 

\maketitle


\section{Introduction}
\label{sec:01}
The present paper is a continuation of the author's work
\cite{yu:gamma}. In loc.~cit. we have determined the number of irreducible
components of a mod $p$ Siegel moduli space with Iwahori level structure. 
The main ingredients are a result of Ng\^o and Genestier
\cite{ngo-genestier:alcoves} that the ordinary
locus is dense in the moduli space, and the surjectivity of a $p$-adic
monodromy due to Faltings and Chai \cite{faltings-chai}, 
also due to Ekedahl \cite{ekedahl:mono}. The goal of this paper
is to investigate the same problem for the non-ordinary locus and 
smaller strata. 


Let $p$ be a rational prime number. Let $N\ge 3$ be a prime-to-$p$
positive integer. We choose a primitive $N$-th root of unity $\zeta_N$ in
$\Qbar\subset \C$ and an embedding $\Qbar\hookrightarrow \Qbar_p$. Let
$\calA_{g,1,N}$ denote the moduli space over $\Z_{(p)}[\zeta_N]$ of
$g$-dimensional 
principally polarized abelian varieties with a full symplectic 
level-$N$ structure with respect to $\zeta_N$. 
The moduli scheme $\calA_{g,1,N}$
has irreducible geometric fibers. Let $\calA$ be the reduction
$\calA_{g,1,N}\otimes \Fpbar$ modulo $p$. For each integer 
$0\le f\le g$, let $\calA^f\subset \calA$ be the locally closed
reduced subscheme that classifies the objects $(A,\lambda,\eta)$ 
whose $p$-rank is $f$. The $p$-rank of an abelian variety $A$ is the
dimension of $A[p](\bar k)$ over 
$\Fp$. It is known due to Koblitz \cite{koblitz:thesis} that each
stratum $\calA^f$ is equi-dimensional of co-dimension $f$ and the
closure of the stratum $\calA^f$ contains $\calA^{f-1}$ for all $f$. 
This result is
generalized to the moduli spaces of arbitrary polarized abelian
varieties by Norman and Oort \cite{norman-oort}


Let $(\calX,\lambda,\eta)\to \calA^f$ be the universal family. 
The maximal etale quotient $\calX[p^\infty]^{\rm et}$ of the 
$p$-divisible group $\calX[p^\infty]$ gives rise
to a $p$-adic monodromy
\[ \rho^f: \pi_1(\calA^f, \bar x)\to \GL_{f}(\Z_p), \]
where $\bar x$ is a geometric point of $\calA^f$. 

In this paper, we prove
\begin{thm}\label{11}
  The homomorphism $\rho^f$ is surjective. 
\end{thm}
The case where $f=g$ is a well-known result proved by Faltings and
Chai \cite{faltings-chai} and independently by Ekedahl
\cite{ekedahl:mono}.
Theorem~\ref{11} answers a question raised in 
Tilouine \cite[Remark below Theorem 2, p. 792]{tilouine:coates}.

A direct consequence of Theorem~\ref{11} is that the associated Igusa
tower over each stratum $\calA^f$ is irreducible except 
when $f=0$ and $g\le 2$ (this is the case where the stratum is 
supersingular, that is, it is entirely contained in the supersingular
locus); see Section~\ref{sec:03}. We apply
Theorem~\ref{11} to the almost ordinary locus of the
moduli spaces with parahoric level structure and determine the number
of irreducible components; see
Section~\ref{sec:04}. In the special case of Iwahori level structure, 
we have the following result:  

Let $\calA_{g,\Gamma_0(p),N}$ denote the moduli space
over $\Z_{(p)}[\zeta_N]$ which parametrizes equivalence classes of
objects $(A,\lambda,\eta,H_\bullet)_S$, where $S$ is a
$\Z_{(p)}[\zeta_N]$-scheme, $(A,\lambda,\eta)$ is in
$\calA_{g,1,N}(S)$, and $H_\bullet$ is a flag of finite 
flat subgroup schemes of $A[p]$
\[ H_\bullet:\quad  0\subset H_1\subset H_2\subset \dots \subset
H_g\subset A[p]\]
such that each $H_i$ is of rank $p^i$ and $H_g$ is isotropic for 
the Weil pairing $e_\lambda$ induced by $\lambda$. 
Let $\calA_{\Gamma_0(p)}:= \calA_{g,\Gamma_0(p),N}\otimes
\Fpbar$ be the reduction modulo $p$ and let
$\calA_{\Gamma_0(p)}^{g-1}:=\calA_{\Gamma_0(p)}\times_
\calA \calA^{g-1}$, the almost ordinary locus of
$\calA_{\Gamma_0(p)}$. We prove (Corollary \ref{43})
\begin{thm}\label{12}
  For $g\ge 2$, the almost ordinary locus $\calA_{\Gamma_0(p)}^{g-1}$
  has $g2^{g-1}$ irreducible components.
\end{thm}


One might expect that the almost ordinary locus
$\calA_{\Gamma_0(p)}^{g-1}$ is dense in the non-ordinary locus
$\calA_{\Gamma_0(p)}^{\rm non-ord}$. If this is true, then it 
would imply the same for the moduli spaces with
{any} parahoric level structure as well (see an argument in
\cite{yu:gamma} for the ordinary case), and then we
could determine the number of irreducible components of these 
non-ordinary loci. However, it is false in general. We examine an
example in Section~\ref{sec:05} (see Proposition~\ref{53}). 

In Section~\ref{sec:045} we show how to use Theorem~\ref{11} to
determine the numbers of  connected components of the $p$-rank strata. 

In Sections \ref{sec:06} and \ref{sec:07}, we give a geometric
characterization of Kottwitz-Rapoport strata for the case $g=2$. The
characterization requires some knowledge on the supersingular
locus. Therefore, a description of the supersingular locus is
included. On the other hand, the characterization also gives us more
information on the supersingular locus through the
Kottwitz-Rapoport stratification. 
This information enables us to determine the number of irreducible 
components of each Kottwitz-Rapoport stratum. 

\section{Proof of Theorem~\ref{11}}
\label{sec:02}

 We may assume that $1\le f<g$ because the case $f=g$ is done in 
 \cite{faltings-chai} and there is nothing to show for $f=0$. Since
 the stratum $\calA^f$ is irreducible (Proposition~\ref{315}), it
 suffices to show the statement for a specific geometric base point. 
 Choose a point $x_0=\ul A_0$ in
 $\calA_{g-f,1,N}\otimes \Fpbar$ whose $p$-rank is zero. Consider the
 morphism
\[ \alpha:\calA_f^{\rm ord} \to \calA^f,\quad \ul A\mapsto \ul A_0\times
\ul A, \]
where $\calA_f^{\rm ord}$ is the ordinary locus of the reduction
$\calA_{f,1,N}\otimes \Fpbar$  mod $p$. Choose a geometric point $\bar
x_1$ of $\calA_f^{\rm ord}$. We have the following commutative diagram for
the $p$-adic monodromies
\[ \xymatrix{
\pi_1(\calA_f^{\rm ord},\bar x_1) \ar[r]^{\alpha_*} \ar[rd]_{\rho^{\rm
    ord}_f} & \pi_1(\calA^f,\bar x_0\times \bar x_1)\ar[d]^{\rho^f} \\
 & \GL_f(\Z_p).
} \]
Since $\rho^{\rm ord}_f$ is surjective \cite{faltings-chai}, $\rho^f$
is also surjective. This completes the proof of Theorem~\ref{11} . 





\section{Irreducibility of the Igusa towers}
\label{sec:03}

\subsection{}

Let $f$ be an integer with $0\le f\le g$. For each integer $m\ge 0$, let
$\calI^f_m$ be the cover of $\calA^f$ over $\Fpbar$ which
parametrizes equivalence classes of objects
$(A,\lambda,\eta,\xi)_S$ where $S$ is an $\Fpbar$-scheme,
$(A,\lambda,\eta)_S$ is in $\calA^f(S)$ and 
$\xi$ is an isomorphism form $\mu_{p^m,S}$ to the multiplicative part
$A[p^m]^{\rm mul}$ of $A[p^m]$ over $S$. 
Let \[ \calI^f:=\{\calI_m^f\}_{m\ge 0} \]
be the Igusa tower over the stratum $\calA^f$. 

\begin{thm}[Oort] \label{32}
   Every non-supersingular Newton polygon stratum of $\calA$ is
   irreducible.  
\end{thm}

This is yet a unpublished result of Oort. See the sketch of his proof
in \cite{oort:bonn}.

\begin{prop}\label{315}
  The stratum $\calA^f$ is irreducible except when $f=0$ and $g\le 2$.
\end{prop}
\begin{proof}
  From all possible symmetric Newton
  polygons, we know that 
  \begin{itemize}
  \item[(a)] the stratum $\calA^f$ is supersingular (means that every
    maximal point of $\calA^f$ is supersingular) if and only if $f=0$
    and $g\le 2$, and
\item [(b)] the stratum $\calA^f$ contains a unique maximal Newton polygon
  stratum as an open dense subset.
  \end{itemize}
Then the proposition follows from Theorem~\ref{32}. \qed  
\end{proof}

\begin{prop}\label{31}
  The Igusa tower $\calI^f$ is irreducible except when $f=0$ and $g\le 2$. 
\end{prop}
\begin{proof}
  The cover $\calI^f_m$ is etale over $\calA^f$ and it represents the
  etale sheaf 
  \[ \ul{Isom}(\mu_{p^m}^{\oplus f}, \calX[p^m]^{\rm mul}),\]
  where $(\calX,\lambda,\eta)\to \calA^f$ is the universal
  family. Therefore, the cover gives rise to the $p$-adic monodromy
  $\rho^f_m: \pi_1(\calA^f,\bar x)\to \GL_f(\Z/p^m \Z)$. 
  By Theorem~\ref{11}, the homomorphism $\rho^f_m$ is surjective for
  all $m\ge 0$. Therefore, each $\calI^f_m$
  is irreducible if the base $\calA^f$ is irreducible. 
  Then the proposition follows from Proposition~\ref{315}. \qed
\end{proof}





When $f=0$, each member $\calI^f_m$ of $\calI^f$ is $\calA^0$; this is
the trivial case. For the non-trivial cases $f\ge 1$, the Igusa towers
are all irreducible.

\subsection{} We consider a variant of the Igusa towers. Let
$f$ be an integer with $1\le f\le g$. For each integer $m\ge 0$, let
$H^f_m:=(\Z/p^m\Z\times \mu_{p^m})^{\oplus f}$ and let $\varphi^f_m:H^f_m\times
H^f_m \to \mu_{p^m}$ be the 
alternating pairing defined by 
\[ \varphi^f_m( (m_i,\zeta_i),(n_i,\eta_i))=\prod_{i=1}^f \eta_i^{m_i}
\zeta_i^{-n_i}, \quad \forall\,m_i,n_i\in \Z/p^m\Z, \ \zeta_i,\eta_i\in
\mu_{p^m}. \] 

Let $\calJ^f_m$ be the cover of $\calA^f$ over $\Fpbar$ which parametrizes
equivalence classes of objects $(A,\lambda,\eta,\xi)_S$ where $S$ is
an $\Fpbar$-scheme, $(A,\lambda,\eta)_S$ is in $\calA^f(S)$ and 
\[ \xi:(H^f_m)_S\to A[p^m] \]
is a monomorphism (both homomorphism and closed
immersion) over $S$ such that 
\[ \varphi^f_m(x,y)=e_\lambda(\xi(x),\xi(y)), \quad \forall\, x,y\in
H^f_m, \]
where $e_\lambda$ is the Weil pairing induced by $\lambda$. Let 
\[ \calJ^f:=\{\calJ^f_m\}_{m\ge 0} \]
be the associated tower over the stratum $\calA^f$.  
\begin{thm}\label{33}
  The tower $\calJ^f$ is irreducible, that is, each member $\calJ^f_m$ is
  irreducible. 
\end{thm}
\begin{proof}
Let $(\calX,\lambda,\iota)\to \calA^f$ be the universal family. 
Consider the canonical filtration  
\[ 0\subset \calX[p^m]^{\rm mul} \subset \calX[p^m]^{0}
\subset \calX[p^m],\]
where $A[p^m]^{0}$ is the neutral connected component of
$A[p^m]$. So we have two canonical short exact sequences
\[ 0\to \calX[p^m]^{\rm mul} \to  \calX[p^m]^{0}\to 
\calX[p^m]^{\rm loc,loc}\to 0,\]
\[ 0\to \calX[p^m]^{0}\to \calX[p^m]\to \calX[p^m]^{\rm et}\to 0.\]
Since these short exact sequences split over a perfect
affine base in \ch $p$, we can find a finite radical surjective morphism
$\pi:\calA'\to A^f$ such that the base change $\calX[p^m]\times_{\calA^f}
\calA'$ admits the canonical decomposition
\begin{equation}
  \label{eq:33}
  \calX[p^m]_{\calA'}=\left(\calX[p^m]_{\calA'}^{\rm mul}\oplus
  \calX[p^m]^{\rm et}_{\rm sub}\right) \oplus  
\calX[p^m]^{\rm loc,loc}_{\rm sub},
\end{equation}
where $\calX[p^m]_{\calA'}^{\rm mul}$ is the base change
$\calX[p^m]_{}^{\rm mul}\times_{\calA^f} \calA'$, 
the middle part $\calX[p^m]^{\rm et}_{\rm sub}$ is the maximal  
etale subgroup scheme of $\calX[p^m]_{\calA'}$ and 
$\calX[p^m]^{\rm loc,loc}_{\rm sub}$ is the maximal 
local-local subgroup scheme of
  $\calX[p^m]_{\calA'}$. Furthermore, we may choose $\calA'$ to be
  irreducible. To see this, let $\calA_0$ be a scheme over
  $\F_q$ such that $\calA^f \simeq \calA_0\otimes_{\F_q}\Fpbar$. Let
  $\{U\}$ be a finite open covering of affine subschemes of
  $\calA_0$. Choose a positive integer $n$ large enough such that
  $\calX[p^m]$ admits the 
  canonical decomposition as (\ref{eq:33}) over $U^{(q^{-n})}$ for
  each $U$, where 
  $F^n:U^{(q^{-n})}\to U$ is the iterated relative Frobenius 
morphism over
  $\F_q$. The subgroup schemes 
\[ \left \{(\calX[p^m]_{U^{(q^{-n})}})^{\rm mul}_{\rm sub}\right
  \}_{U^{(q^{-n})}}, \quad \text{(resp. 
$\left \{(\calX[p^m]_{U^{(q^{-n})}})^{\rm loc,loc}_{\rm sub}\right
\}_{U^{(q^{-n})}}$ )} \] 
glue to a subgroup scheme 
$\calX[p^m]^{\rm mul}_{\rm sub}$
(resp. $\calX[p^m]^{\rm loc,loc}_{\rm sub}$) over
$\calA_0^{(q^{-n})}$. Clearly $\calA_0^{(q^{-n})}$ is irreducible. 
We may take $\calA':=\calA_0^{(q^{-n})}\otimes_{\F_q}{\Fpbar}$ 
and let $\pi=F^n$, then $\calX[p^m]_{\calA'}$ 
admits the canonical decomposition.

Let $\calJ'_m$ be the etale cover of $\calA'$ that represents the
etale sheaf
\[ \calP_m:=\ul {Isom}\left ( (H^f_m,\varphi^f_m)_{\calA'},
(\calX[p^m]_{\calA'}^{\rm mul}\oplus \calX[p^m]^{\rm et}_{\rm sub},
e_\lambda)\right ). \]
Since any section $\xi$ of $\calP_m$ is determined by its restriction
$\xi$ on $(\Z/p^m)^{\oplus f}$, the restriction map 
$\xi\mapsto \xi|_{(\Z/p^m)^{\oplus f}}$ gives an 
isomorphism 
\[ \calP_m\simeq \ul{Isom}((\Z/p^m)^{\oplus f}_{\calA'},\calX[p^m]^{\rm
  et}_{\rm sub}). \]
Therefore, $\calJ'_m$ corresponds to the $p$-adic monodromy
$\rho'_m: \pi_1(\calA',\bar x')\to\GL_f(\Z/p^m\Z)$ and we have the
commutative diagram
\[ \xymatrix{
\pi_1(\calA',\bar x') \ar[r]^{\pi_*} \ar[rd]_{\rho'_m} & 
\pi_1(\calA^f,\bar x)\ar[d]^{\rho^f_m} \\
 & \GL_f(\Z/p^m\Z).
} \] 
Since $\pi_*$ is an isomorphism and $\rho^f_m$ is surjective, 
$\rho'_m$ is surjective. Therefore, $\calJ'_m$ is 
irreducible as $\calA'$ is so. Let $(\calX',\lambda',\eta')$ 
be the base change of $(\calX,\lambda,\eta)$ 
over $\calA'$ and let $\xi'\in
\calP_m(\calJ'_m)$ be the universal section. Then the family
$(\calX',\lambda',\eta',\xi')\to\calJ'_m$ gives rise to a morphism
$\alpha:\calJ'_m \to \calJ^f_m$. Clearly this map is surjective, hence
$\calJ^f_m$ is irreducible. \qed 
\end{proof}
\begin{remark} 
The irreducibility of $\calJ_1^g$  is studied in
  \cite{yu:gamma} (that was denoted $\calA^{\rm
    ord}_{\Gamma(p)}$ there). The lines 1-2 of p.~2593 in loc.~cit.~are
  incorrect. The moduli scheme  
$\calA^{\rm  ord}_{\Gamma(p)}$ is not etale over $\calA^{\rm
  ord}_{g,1,N}$ because the extension
\[ 0\to \calX[p]^0\to \calX[p]\to \calX[p]^{\text{\rm et}}\to 0 \]
does not split over any finite etale base change. However, this does
not effect the conclusion on irreducibility of $\calA^{\rm
  ord}_{\Gamma(p)}$; we just need a modified argument as in the proof
above. 
\end{remark}
\section{The almost ordinary locus of the moduli spaces with parahoric
level structure}
\label{sec:04}

\subsection{}
\label{sec:41}
We keep the notation as before.
Let $\ul k=(k_1,\dots,k_r)$ be a tuple of positive integers $k_i\ge 1$ with
$\sum_{i=1}^r k_i \le g$. Set $h(i):=\sum_{j=1}^i k_j$ for
$1\le i\le r$ and $h(0)=0$. Let $\calA_{g,\ul k,N}$ denote the moduli
space over $\Z_{(p)}[\zeta_N]$ that parametrizes equivalence classes 
of objects
$(A,\lambda,\eta, H_\bullet)_S$, where $S$ is a
$\Z_{(p)}[\zeta_N]$-scheme, $(A,\lambda,\eta)$ is in
$\calA_{g,1,N}(S)$, and $H_\bullet$ is a flag of finite subgroup
schemes of $A[p]$ 
\[ H_\bullet\quad 0=H_{h(0)}\subset H_{h(1)}\subset \dots\subset
H_{h(r)}\subset A[p] \]
such that $H_{h(i)}$ is locally free of rank $p^{h(i)}$ and $H_{h(r)}$ is
isotropic for the Weil pairing $e_{\lambda}$ induced by the
polarization $\lambda$. When $r=g$, the moduli scheme
$\calA_{g,\ul k,N}$ is $\calA_{g,\Gamma_0(p),N}$
defined in Section 1. 

Let $\calA_{\ul k}:=
\calA_{g,{\ul k},N} \otimes \ol \F_p$ be the reduction
modulo $p$. For $0\le f\le g$, let 
$\calA_{\ul k}^{f}:=\calA_{\ul  k}\times_\calA \calA^{f}$, 
the $p$-rank $f$ stratum of
the moduli space $\calA_{\ul k}$. 
For an $\Fpbar$-scheme $S$, the
$S$-valued set $\calA_{\ul k}^{f}(S)$ consists of objects
$(A,\lambda,\eta,H_\bullet)_S$ in $\calA_{\ul  k}(S)$ such that the
canonical morphism $S\to \calA$ given by the family
$(A,\lambda,\eta)_S$ factors through the subscheme $\calA^f$.
Note that from the definition one can not determine
whether $\calA_{\ul k}^{f}$ is reduced.  
We will compute the number of irreducible
components of $\calA_{\ul k}^{f}$ in the case where $f=g-1$, the
almost ordinary locus.

We first seek discrete invariants for geometric points on 
$\calA_{\ul  k}^{g-1}$. Let $k$ be an \ac field of \ch $p$. Fix 
a supersingular elliptic curve $E_0$ over $\Fpbar$. Let 
$(A,\lambda,\eta, H_\bullet)$ be a point of $\calA^{g-1}_{\ul k}(k)$. We
  have
\[ A[p]\simeq (\Z/p\Z \times \mu_p)^{g-1}\times E_0[p]. \]
There are two cases:
\begin{itemize}
\item [(a)] $H_{h(r)}$ does not have local-local part. This occurs
  only when $h(r)<g$. For each $1\le i\le r$, the finite group scheme
  $H_h(i)/H_{h(i-1)}$ has the form $\mu_p^{\tau(i)}\times
  (\Z/p\Z)^{k_i-\tau(i)}$ for a non-negative integer $0\le \tau(i)\le
    k_i$.  
\item [(b)] $H_{h(r)}$ has non-trivial local-local part. 
  There is a unique integer $1\le j \le r$ such that $H_{h(j)}$ 
  has non-trivial local-local part 
  and $H_{h(j-1)}$ has no local-local part. For each $1\le i\le
  r$, the finite group scheme 
  $H_{h(i)}/H_{h(i-1)}$ has the form 
\[ 
\begin{cases}
  \mu_p^{\tau(i)}\times (\Z/p\Z)^{k_i-\tau(i)}\ \text{for some integer $0\le 
  \tau(i)\le k_i$} & \text{if $i\neq j$};\\
  \mu_p^{\tau(i)}\times (\Z/p\Z)^{k_i-1-\tau(i)}\times \alpha_p\  
   \text{for some integer $0\le \tau(i)\le 
    k_i-1$} & \text{if $i=j$}.
\end{cases}
\]
\end{itemize}

\subsection{}
\label{sec:42}
For each  $m\ge 0$, define $I(m):=[0,m]\cap\Z$. Set 
\[ I^0_{\ul k}:=
\begin{cases}
  \emptyset & \text{if $h(r)=g$};\\
   \{0\}\times\prod_{i=1}^r I(k_i) & \text{if $h(r)<g$}. 
\end{cases} \]
For each $1\le j\le r$, set 
\[ I^j_{\ul k}:=\{\ul\tau=(j, \tau(1),\dots,\tau(r)); \quad \tau(i)\in
I(k_i) \text{ for $i\neq j$ and } \tau(j)\in I(k_j-1)\,\}. \]
The finite set $I^0_{\ul k}$ (resp. $I^j_{\ul k}$  for $j>0$) will be used to
parameterize discrete invariants in the case (a) (resp. the case (b)).

For $\ul \tau=(0, \tau(1),\dots,\tau(r))\in I^0_{\ul k}$, we say a
geometric point $\ul A$ in $\calA_{\ul k}^{g-1}$ is of type $\ul \tau$ if 
\begin{itemize}
\item $H_{h(r)}$ has no local-local part, and
\item the multiplicative part of $H_{h(i)}/H_{h(i-1)}$ is of rank
  $p^{\tau(i)}$ for all $1\le i\le r$. 
\end{itemize}
For $\ul \tau=(j, \tau(1),\dots,\tau(r))\in I^j_{\ul k}$, where $1\le
j\le r$, 
we say a geometric point $\ul A$ in $\calA_{\ul k}^{g-1}$ is of 
type $\ul \tau$ if \begin{itemize}
\item $H_{h(j)}$ has local-local part and $H_{h(j-1)}$ has no
local-local part, and
\item the multiplicative part of $H_{h(i)}/H_{h(i-1)}$ is of rank
  $p^{\tau(i)}$ for all $1\le i\le r$. 
\end{itemize}
Then we have an assignment $\ul A \mapsto \ul \tau(\ul A)$, 
which gives a surjective map 
\[ \tau: \calA^{g-1}_{\ul k} \to \coprod_{0\le j\le r} I^j_{\ul k}. \]
It is easy to see that this map is locally constant for the Zariski
topology (use the argument in the proof of Theorem~\ref{33}).  
For a fixed type $\ul \tau$, let $\calA^{g-1}_{\ul k,\ul \tau}$ be the
union of the connected components of $\calA^{g-1}_{\ul k}$ 
whose objects are of type $\ul \tau$. 
We write $\calA^{g-1}_{\ul k}$ into a disjoint
union of open subschemes  
\[ \calA^{g-1}_{\ul k}=\coprod_{0\le j \le r}\ \coprod_{\ul \tau\in
  I^j_{\ul k}} \calA^{g-1}_{\ul k,\ul \tau}.
\]
\begin{thm}\label{41}
  For each integer $0\le j\le r$ and each type $\ul \tau\in I^j_{\ul k}$,
  there is a finite surjective morphism $\pi_{\ul \tau}: \calJ^{g-1}_1\to
  \calA^{g-1}_{\ul k,\ul \tau}$. Consequently, each stratum
  $\calA^{g-1}_{\ul k,\ul \tau}$ is irreducible
  for $g\ge 2$. 
\end{thm}
\begin{proof}
  Let $(\calX,\lambda,\eta,\xi)\to \calJ^{g-1}_1$ be the universal
  family. The image $\xi(H^{g-1}_1)$ is the etale-multiplicative part
  $\calX[p]^{\rm em}$
  of $\calX[p]$, and the orthogonal complement of $\xi(H^{g-1}_1)$ for
  the Weil pairing $e_\lambda$ is the local-local part $\calX[p]^{\rm
  loc,loc}$ of $\calX[p]$. Namely, we have 
\[ \calX[p]=\calX[p]^{\rm em} \times \calX[p]^{\rm loc,loc},\quad
  \xi:(\mu_p\times \Z/p\Z)^{g-1}\stackrel{\sim}{\to} \calX[p]^{\rm
  em}. \]
Let $C$ be the kernel of the relative Frobenius morphism
  \[ F_{\calX/\calJ^{g-1}_1}:\calX[p]^{\rm loc,loc}\to (\calX[p]^{\rm
  loc,loc})^{(p)}. \]
\begin{itemize}
\item [(a)] If $\tau\in I^0_{\ul k}$, let $K_i:=\mu_p^{\tau(i)}\times
  (\Z/p\Z)^{k_i-\tau(i)}$ for $1\le i\le r$. For $1\le m\le r$, set
    $H_{h(m)}:=\xi(\prod_{i=1}^m K_i)$. Then we define a family
    $(\calX,\lambda,\eta,H_\bullet)\to \calJ^{g-1}_1$ and this family
  induces a natural morphism $\pi_{\ul \tau}:\calJ^{g-1}_1\to \calA_{\ul k}$. 
\item [(b)] If $\tau\in I^j_{\ul k}$ for some $1\le j\le r$, let
\[ K_i:=
\begin{cases}
  \mu_p^{\tau(i)}\times
  (\Z/p\Z)^{k_i-\tau(i)} & \text{if $i\neq j$}; \\
\mu_p^{\tau(i)}\times
  (\Z/p\Z)^{k_i-1-\tau(i)} & \text{if $i= j$}.  
\end{cases}\]
For $1\le m\le r$, set
\[ H_{h(m)}:=
\begin{cases}
  \xi(\prod_{i=1}^m K_i) & \text{if $m<j$}; \\
\xi(\prod_{i=1}^m K_i)\times C  & \text{if $m\ge j$}.  
\end{cases}\]
Then we define a family
    $(\calX,\lambda,\eta,H_\bullet)\to \calJ^{g-1}_1$ and this family
    induces a natural
    morphism $\pi_{\ul \tau}:\calJ^{g-1}_1\to \calA_{\ul k}$.
\end{itemize}
It is clear that $\pi_{\ul \tau}$ factors through the almost ordinary locus
$\calA^{g-1}_{\ul k}$. Moreover, the image lands in the open subscheme
$\calA^{g-1}_{\ul k,\ul \tau}$ by the construction. So we get the
morphism $\pi_{\ul \tau}:\calJ^{g-1}_1\to\calA^{g-1}_{\ul k,\ul \tau}$. One checks
easily that $\pi_{\ul \tau}$ is surjective. Since the composition
$\calJ^{g-1}_1\to\calA^{g-1}_{\ul k,\ul \tau}\to \calA^{g-1}$ is finite,
$\pi_{\ul \tau}$ is finite. This completes the proof. \qed
\end{proof}

Note that in the proof we use the universal family to define the
morphism $\pi_{\ul \tau}$ instead of defining
$\pi_{\ul \tau}(x)$ pointwisely. The reason is that $\calJ^{g-1}_1$ or
$\calA^{g-1}_{\ul k}$ (defined by the fiber product) could be
non-reduced. 

\begin{cor}\label{42}
  For $g\ge 2$, the almost ordinary stratum $\calA^{g-1}_{\ul k}$ has 
\[ \sum_{j=0}^r |I^j_{\ul k}|=(k_1+1)\dots(k_r+1)\left
  [\epsilon+\frac{k_1}{k_1+1}+\dots+\frac{k_r}{k_r+1}\right] \] 
irreducible components, where $\epsilon=1$ if $h(r)<g$ and $\epsilon=0$
if $h(r)=g$.
\end{cor} 

For the Iwahori case, $r=g$ and $k_i=1$ for all $i$, so
Corollary~\ref{42} gives
 
\begin{cor}\label{43}
  For $g\ge 2$, the almost ordinary locus $\calA_{\Gamma_0(p)}^{g-1}$
  has $g2^{g-1}$ irreducible components.
\end{cor}

\section{Connected components of $p$-rank strata}
\label{sec:045}

In the previous section we study irreducible components of the almost
ordinary locus. In this section we consider lower $p$-rank strata. 
We know that when $g\ge 2$ and $0\le f\le g-2$, the natural morphism
$\calA_{\ul k}^f \to \calA^f$ is not finite in general. This limits
the method of using $p$-adic monodromy to the irreducibility problem
in the present case. The obstacle results from the
fibration ``moving $\alpha_p$-subgroups''. If one contracts the
fibration, then one obtains a finite morphism for which the $p$-adic
monodromy results can be applied. Proceeding this approach, we
obtain information on connected components instead.  \\

Keep the notation in the previous sections. Assume that $g\ge 2$.

\subsection{}\label{sec:0451}
Fix a tuple $\ul k=(k_1,\dots,k_r)$ of positive integers
with $\sum_{i=1}^r k_i\le g$ and an integer $f$ with $0\le f\le g$ 
as before. Again, we first seek discrete invariants for 
geometric points in $\calA^f_{\ul k}$. 
Let $(A,\lambda,\eta, H_\bullet)$ be a point in $\calA^f_{\ul k}(k)$. 
We have a decomposition 
\[ H_{h(i)}=\left ( H^{\rm et}_{h(i)}\oplus H^{\rm mul}_{h(i)} 
  \right ) \oplus H^{\rm loc,loc}_{h(i)} \] 
into etale-multiplicative part and local-local part. 
Suppose that $H^{\rm et}_{h(i)}\oplus H^{\rm mul}_{h(i)}$ has rank
$p^{a(i)}$ and $H^{\rm loc,loc}_{h(i)}$ has rank $p^{b(i)}$. 
Put $a(0)=0$ and $m(i):=a(i)-a(i-1)$ for each $1\le i\le r$. 
It is easy to see that
\begin{equation}
  \label{eq:451}
    0\le m(i)\le k_i, \quad \text{and \ } 
   f-(g-h(r))\le \sum_{i=1}^r m(i) \le f,  
\end{equation}
where $h(i)=\sum_{j=1}^i k_j$ as before. 
Let $G_i:=H^{\rm et}_{h(i)}\oplus H^{\rm mul}_{h(i)}$. Then 
the successive quotient $G_i/G_{i-1}$ has rank $p^{m(i)}$. 
Let 
\[ \tau(i):=\log_p \rank (G_i/G_{i-1})^{\rm mul}, \quad \forall\, 1\le i\le
r. \]
We have $0\le \tau(i)\le m(i)$ for all $i$. We call the pair of
$r$-tuples 
\[ (\ul m, \ul \tau)=[(m(1),\dots, m(r)),(\tau(1), \dots
\tau(r))] \] 
the {\it graded etale-multiplicative type} associated to the
object $(A,\lambda,\eta, H_\bullet)$, abbreviated as {\it gem} type.

Conversely, fix a tuple of integers $(m(1),\dots, m(r))$ satisfying
(\ref{eq:451}). Let $(A,\lambda,\eta)$ be an object in $\calA^f(k)$
and $G_\bullet$ a flag of finite flat subgroup schemes  
\[ 0=G_0\subset G_1\subset\dots\subset
G_r\subset A[p] \] 
such that 
\begin{itemize}
\item [(1)] the group scheme $G_r$ is isotropic with respect to the
Weil pairing $e_\lambda$, and
\item [(2)] each $G_i$ has no local-local part and the 
quotient $G_i/G_{i-1}$ has rank $p^{m(i)}$.
\end{itemize}
Then one can lift to an
object $(A,\lambda, \eta, H_\bullet)\in \calA^f_{\ul k}(k)$ such that
$H^{\rm et}_{h(i)}\oplus H^{\rm mul}_{h(i)}=G_i$ for all $i$.

\subsection{}
\label{sec:0452}
Define 
\begin{equation}
  \label{eq:452}
  \begin{split}
    \Sigma_0(\ul k,f):=  \{\, & \ul m=(m(1),
     \dots, m(r))\in \Z^r\, |  \\
    &  0\le m(i)\le k_i, \forall\, i, \text{ and } 
     f-(g-h(r))\le \sum_i^r m(i) \le f \, \}, 
  \end{split}
\end{equation}
\begin{equation}
  \label{eq:453}
  \begin{split}
\Sigma(\ul k,f):=\{\, & (\ul m,\ul \tau)=((m(1),\dots, m(r)),
(\tau(1),\dots, \tau(r))) \in \Z^r\times\Z^r \, | \\
& \, \ul m \in \Sigma_0(\ul k,f), \text{ and } 0\le \tau (i)\le m(i),\ 
\forall\, i \, \}. 
\end{split}
\end{equation}

We have a natural map from $\Sigma(\ul k,f)$ to $\Sigma_0(\ul k,f)$
sending $(\ul m,\ul \tau)$ to $\ul m$. The finite set $\Sigma(\ul
k,f)$ is exactly that of all possible graded etale-multiplicative
types of points in $\calA_{\ul k}^f$. 

For any element $\ul m \in \Sigma_0(\ul k,f)$, we define a scheme
$T(\ul m)$ over $\Fpbar$ as follows. For any locally Noetherian
$\Fpbar$-scheme $S$, 
the $S$-valued set $T(\ul m)(S)$ classifies equivalence classes of
objects $(A,\lambda,\eta,G_\bullet)_S$, where 
\begin{itemize}
\item $(A,\lambda,\eta)_S$ is in $\calA^f(S)$, and 
\item $G_\bullet$ is a flag of finite flat subgroup schemes 
\[ G_\bullet:\quad 0=G_0\subset G_1\subset\dots\subset G_r\subset 
A[p]  \]
\end{itemize}
satisfying the conditions (1) and (2) above. 

For any element $(\ul m,\ul \tau)\in \Sigma(\ul k,f)$, let $T(\ul m,
\ul \tau)$ be the (open) subscheme of $T(\ul m)$ that 
consists of objects
$(A,\lambda,\eta,G_\bullet)_S$ such that the multiplicative part of
the quotient $G_{i}/G_{i-1}$ has rank $p^{\tau(i)}$ for all $1\le i\le
r$. Put $T^f_{\ul k}:=\coprod_{\ul m\in \Sigma_0(\ul k,f)} T(\ul m)$. 
Clearly, we have 
\[ T^f_{\ul k}=\coprod_{(\ul m,\ul \tau)\in\Sigma(\ul k,f)} T(\ul
m,\ul \tau). \]

Let $(\calX,\lambda, \eta,\wt H_\bullet)\to 
\calA^f_{\ul k}$ be the universal family over $\calA^f_{\ul k}$. 
Then there is a finite 
dominant homeomorphic morphism
$\pi: \calA'\to \calA^f_{\ul k}$ such that the base change 
$\calX[p]_{\calA'}$ admits the canonical decomposition
\begin{equation}
  \label{eq:454}
  \calX[p]_{\calA'}=(\calX[p]^{\rm mul}\oplus \calX[p]^{\rm et})\oplus
  \calX[p]^{\rm loc,loc}
\end{equation}
into etale-multiplicative part and local-local part; see the proof of
Theorem~\ref{33}. Accordingly, we have the same decomposition
 \begin{equation}
  \label{eq:455}
  \wt H_{h(i),\calA'}=(\wt H_{h(i)}^{\rm mul}\oplus 
   \wt H_{h(i)}^{\rm et})\oplus
  \wt H_{h(i)}^{\rm loc,loc}
\end{equation}
for all $i$.

Since these subgroup schemes are locally free, their ranks are
constant on each connected component of $\calA'$. 
Let $\calA'(\ul m,\ul \tau)\subset \calA'$ be the union of the
connected components whose objects are of gem type 
$(\ul m,\ul \tau)$. 
Let $\calA(\ul m, \ul \tau)\subset \calA^f_{\ul k}$ be 
the open subscheme 
$\pi(\calA'(\ul m,\ul \tau))$. In particular, $\pi: \calA'(\ul m,\ul
\tau)\to \calA(\ul m, \ul \tau)$ is a homeomorphism.
Again, we have
\begin{equation}
  \label{eq:456}
  \calA^f_{\ul k}=\coprod_{(\ul m, \ul \tau)} \calA(\ul m, \ul \tau), 
  \quad \calA'=\coprod_{(\ul m, \ul \tau)} \calA'(\ul m, \ul \tau),
\end{equation}
where $(\ul m,\ul \tau)$ runs through all elements in 
$\Sigma(\ul k,f)$. 

Consider the universal family $(\calX, \lambda, \eta, \wt
H_\bullet)_{\calA'(\ul m,\ul \tau)}$ restricted on the open 
subscheme $\calA'(\ul m,\ul \tau)$. Put 
\[ \wt G_i:=\wt H_{h(i)}^{\rm mul}\oplus 
   \wt H_{h(i)}^{\rm et}, \] 
for all $i$. We get a family 
$(\calX, \lambda, \eta, \wt G_\bullet)_{\calA'(\ul m,\ul \tau)}$. 
This gives rise to a natural morphism 
\begin{equation}
  \label{eq:457}
  c(\ul m,\ul \tau): \calA'(\ul m,\ul \tau)\to T(\ul m,\ul \tau),
\end{equation}
which is proper and surjective (see Subsection~\ref{sec:0451}). 

If $f=0$, then the set $\Sigma(\ul k,f)$ consists of only one element
$(\ul 0,\ul 0)$. In this case, we have 
\begin{equation}
  \label{eq:458}
  \calA'=\calA^0_{\ul k}, \quad T^0_{\ul k}=\calA^0,\quad \text{and}
  \quad c(\ul 0, \ul 0):\calA^0_{\ul k}\to \calA^0. 
\end{equation}

\begin{prop}\label{451} \

\npr {\rm (1)} The stratum $\calA^0$ is connected and it is
irreducible if $g\ge 3$.

\npr {\rm (2)} Suppose that $f\ge 1$. Then there is a finite
surjective morphism $\calJ^f_1 \to T(\ul m,\ul \tau)$. Consequently,
each scheme $T(\ul m,\ul \tau)$ is irreducible. 
\end{prop}

\begin{proof}
  (1) When $g=2$, this is a special case of Proposition 7.3 in Oort
      \cite{oort:eo}. When 
      $g\ge 3$, this is obtained in the proof of Proposition~\ref{31}
      using Theorem~\ref{32}, a result of Oort. 

  (2) The construction is similar to that as in
     Theorem~\ref{41}. Therefore, we do not repeat it. \qed
\end{proof}

\subsection{}
\label{sec:0453}
Let $\ul n=(n(1), \dots n(r'))$ be a tuple of positive integers with
$\sum_{i=1}^{r'} n(i)\le g-f$. We may identify 
the moduli space $\calA_{\ul n}$ with the moduli space 
that parametrizes equivalence classes of chains of isogenies 
\[ \begin{CD}
  \ul A_\bullet: \quad \ul A_0 @>\alpha_{1}>> \ul A_1 @>>> 
  \dots @>>>\ul A_{r'-1} @>\alpha_{r'}>>\ul A_{r'},
\end{CD}\]
where 
\begin{itemize}
\item each $\ul A_i=(A_i,\lambda_i,\eta_i)$ is 
a polarized abelian scheme with a symplectic level-$N$ structure,
\item $\ul A_0$ is an object in $\calA$, and 
\item each $\alpha_i$ is an isogeny of degree $p^{n(i)}$ that
  preserves the level 
structures and the polarizations except when $i=1$, and in this case
one has $\alpha_1^*\lambda_1=p\lambda_0$. 
\end{itemize}

Define $W_{\ul n}\subset \calA_{\ul n}$ to be the reduced subscheme
consisting of objects $\ul A_\bullet$ such that the kernel $\ker
\alpha_i$ is of local-local type for all $i$. 

Let $(\ul m,\ul \tau)$ be an element in $\Sigma(\ul k,f)$.
Let $a(i):=\sum_{j=1}^i m(j)$ and $b(i):=h(i)-a(i)$. Put
$a(0)=b(0)=0$. Write the set $\{b(i); 0\le i\le r\}$ as
\[ \{b'(0), b'(1),\dots b'(r')\} \quad \text{with}\
0=b'(0)<b'(1)<\dots <b'(r'). \]
We have $r'\le r$, $r'\le b'(r')$ and $h(r)-f\le b'(r')\le g-f$. Set
$n(i):=b'(i)-b'(i-1)$ for $1\le i\le r'$. So we define a tuple $\ul n$
of integers from the pair $(\ul m,\ul \tau)$, and have a scheme
$W_{\ul n}$.

Let $x=(A_x,\lambda_x,\eta_x, G_\bullet)$ be a point in  
$T(\ul m,\ul \tau)(k)$. The fiber $c(\ul m,\ul \tau)^{-1}(x)$ 
consists of flags of finite flat subgroup schemes
\[ 0=K_{b(0)}\subset K_{b(1)}\subset\dots\subset 
K_{b(r)}\subset A_x[p] \]
such that $K_{b(r)}$ is isotropic for the Weil pairing $e_{\lambda_x}$
and each $K_{b(i)}$ is local-local of rank $p^{b(i)}$. This is the
same as flags of finite flat subgroup schemes
\[ 0=K_{b'(0)}\subset K_{b'(1)}\subset\dots 
\subset K_{b'(r')}\subset A_x[p] \] 
with the same properties. This proves $c(\ul m,\ul
\tau)^{-1}(x)=W_{\ul n}(x)$, where 
\[ W_{\ul n}(x):=\{\ul A_\bullet \in W_{\ul n}\, ;\, \ul
A_0=(A_x,\lambda_x,\eta_x)\, \}. \]

It is easy to see that the reduced scheme $W_{\ul n}(x)$ is 
connected. Indeed, for $1\le d \le g-f$, let $\ul 1^d:=(1,\dots, 1)$
with length $d$. We see that the fiber of the natural morphism 
$W_{\ul 1^i}(x)\to W_{\ul 1^{i-1}}(x)$ is a projective space, as it is
the family of $\alpha_p$-subgroups in $\ker \lambda_{i-1}$. This shows
that the scheme $W_{\ul 1^{g-f}}(x)$ is connected. 
Since the forgetful morphism $W_{\ul 1^{g-f}}(x)\to W_{\ul n}(x)$ is
  surjective, the scheme $W_{\ul n}(x)$ is connected. In
  conclusion, we have proved 

\begin{prop}\label{452}
  Any fiber of the morphism $c(\ul m,\ul \tau)$ (\ref{eq:457}) is connected.  
\end{prop}

\begin{thm}\label{453}
  Every open subscheme $\calA(\ul m, \ul \tau)\subset \calA^f_{\ul k}$
  is connected. Consequently, the $p$-rank $f$ stratum $\calA^f_{\ul
  k}$ has $|\Sigma(\ul k,f)|$ connected components.
\end{thm}
\begin{proof}
  It follows from Propositions~\ref{451} and~\ref{452} that
  $\calA'(\ul m,\ul \tau)$ is connected. Since  $\calA(\ul m,\ul
  \tau)$ is homeomorphic to $\calA'(\ul m,\ul \tau)$, it is 
  connected. The second statement follows from (\ref{eq:456}). \qed
\end{proof}

\begin{remark}
The scheme $T^f_{\ul k}$ is closely related to the Stein factorization
of the natural morphism $\calA^f_{\ul k}\to \calA^f$. Indeed, let $T$
(resp. $T'$) be 
the Stein factorization of $\calA^f_{\ul k}\to \calA^f$ (resp. of
$\calA'\to \calA^f$). Then there are natural finite morphisms
$\pi_1:T'\to T$ and $\pi_2:T'\to T^f_{\ul k}$. It is not hard to see
that  these morphisms are homeomorphic. 
In some sense, $T^f_{\ul k}$ provides a 
``modular interpretation'' of the scheme $T$.   
\end{remark}




\section{The Siegel 3-fold with Iwahori level structure}
\label{sec:05}

In this section we describe the Kottwitz-Rapoport stratification on
the Siegel $3$-fold $\calA_{2,\Gamma_0(p)}$ with Iwahori level
structure. Our references are de Jong \cite{dejong:gamma}, Kottwitz
and Rapoport \cite{kottwitz-rapoport:alcoves}, 
T.~Haines \cite{haines:clay}, and Ng\^o and Genestier
\cite{ngo-genestier:alcoves}. The geometric part of Tilouine 
\cite{tilouine:coates} is also helpful to us. 
Then we conclude the following results as consequences:
\begin{itemize}
\item [(a)] The almost ordinary locus
  $\calA_{2,\Gamma_0(p)}^{1}$ is not dense in the non-ordinary locus
  $\calA_{2,\Gamma_0(p)}^{\rm non-ord}$.
\item [(b)] The supersingular locus
  $\calS_{2,\Gamma_0(p)}$ of
  $\calA_{2,\Gamma_0(p)}$ is not equi-dimensional. It consists of
  both one-dimensional components and two-dimensional components.
\end{itemize}

\subsection{Local models}
\label{sec:51}
The Kottwitz-Rapoport stratification is defined through local
model diagrams. We setup a lattice chain which defines a local model.
Let $\calO$ be a complete discrete valuation ring, $K$ its fraction
field, $\pi$ an uniformizer of $\calO$, and $\kappa:=\calO/\pi\calO$ the
residue field. We require that $\mathrm{char} \kappa   = p>0$. Set
$V:=K^{2n}$ and let 
$e_1,\dots,e_{2n}$ be the standard basis. Denote by $\psi:V\times V\to
K$ the non-degenerate alternating form whose non-zero pairings are 
\[ \psi(e_i,e_{2n+1-i})=1, \quad 1\le i\le n, \]
\[ \psi(e_i,e_{2n+1-i})=-1, \quad i\ge n+1,\]
The representing matrix for
$\psi$ is 
\[ 
\begin{pmatrix}
  0 & \wt I \\
  -\wt I & 0 
\end{pmatrix},\quad \wt I=\text{anti-diag(1,\dots,1)}.\]
Let $\GSp_{2n}$ be the reductive algebraic group of symplectic
similitudes with respect to $\psi$. 
Let 
\[ \pi L_0=L_{-2n}\subset L_{-2n+1}\subset \dots\subset L_{-1}\subset
L_0=\calO^{2n} \]
be a chain of $\calO$-lattices in $V$ where the lattice $L_{-i}$ is
generated 
by $e_1,\dots, e_{2n-i}$, 
$\pi e_{2n-i+1},\dots,\pi e_{2n}$. The
$\calO$-submodule in $V$ generated by $x_1,\dots, x_k\in V$ is denoted
by $<x_1,\dots, x_k>$. 
Thus $L_{i-2n}=<e_1,\dots 
, e_{i},\pi e_{i+1},\dots,\pi e_{2n}>$. 

For $0\le i\le 2n$, let $\Lambda_{i-2n}=\calO^{2n}$ and define
$\beta_{i-2n}:\Lambda_{i-2n}\to \Lambda_{i-2n+1}$ for $i<2n$ by 
\[ \beta_{i-2n}(e_{i+1})=\pi e_{i+1},\quad \text{and}\quad
\beta_{i-2n}(e_{j})= e_{j}\quad \text{for $j\neq i+1$}. \]  
We have 
\[ 
\begin{CD}
  \Lambda_{-2n}@>\beta_{-2n}>> \Lambda_{-2n+1}@>>> \dots @>>>\Lambda_{-1}
  @>\beta_{-1}>>\Lambda_0,
\end{CD}\]
and there is a unique isomorphism $a_{-i}:\Lambda_{-i}\to L_{-i}$ with
$a_0={\rm id}$ such that
the diagram 
\[ 
\begin{CD}
  \Lambda_{-i}@>\beta_{-i}>> \Lambda_{-i+1}\\
@VVa_{-i}V  @VVa_{-i+1}V \\ 
L_{-i}@>{\rm incl}>> L_{-i+1} 
\end{CD}\] 
commutes.
Let $\wt \beta_{-i}:\Lambda_{-i}\to \Lambda_0$ denote the composition of
the morphisms
\[ 
\begin{CD}
  \Lambda_{-i}@>\beta_{-i}>> \Lambda_{-i+1}@>>> \dots @>>>\Lambda_{-1}
  @>\beta_{-1}>>\Lambda_0.
\end{CD}\]
Let $\psi_0=\psi$ be on the form on $\Lambda_0=L_0$. There is a
perfect non-degenerate alternating form $\psi_{-n}$ on $\Lambda_{-n}$
such that 
\[ \psi_0(\wt \beta_{-n}(x), \wt \beta_{-n}(y))=\pi\psi_n(x,y), \quad
\forall\, x, y\in \Lambda_{-n}. \]

Let ${\bfM}^{\rm loc}$ denote the projective $\calO$-scheme that
represents the functor which sends an $\calO$-scheme $S$ to the
set of the collections of locally free $\calO_S$-submodules
$\scrF_{-i}\subset \Lambda_{-i}\otimes \calO_S$ of rank $n$ 
for $0\le i\le n$ such that 
\begin{itemize}
\item [(i)] $\scrF_0$ and $\scrF_{-n}$ are isotropic with respect to
  $\psi_0$ and $\psi_{-n}$, respectively.
\item [(ii)] $\scrF_{-i}$ locally is a direct summand of
  $\Lambda_{-i}\otimes \calO_S$ for all $i$,
\item [(iii)] $\beta_{-i}(\scrF_{-i})\subset \scrF_{-i+1}$ for all $i$.
\end{itemize}

By an automorphism on $\Lambda_\bullet\otimes \calO_S$, where $S$ is
an $\calO$-scheme, we mean a collection of
automorphisms $g_{-i}$ on $\Lambda_{-i}\otimes \calO_S$ such that
$g_{-i}$ commutes with the morphisms $\beta_{-i}$ for all $i$ and
$g_0$ and $g_{-n}$ preserve the forms $\psi_0$ and $\psi_{-n}$, 
respectively, up to invertible scalars. We denote by
$\mathrm{Aut}(\Lambda_\bullet\otimes \calO_S, \psi_0,\psi_{-n})$ the
group of automorphisms 
on $\Lambda_\bullet\otimes \calO_S$. 

Let $\calG$ be the group
scheme over $\calO$ that represents the functor
\[ S\mapsto \mathrm{Aut}(\Lambda_\bullet\otimes \calO_S, 
\psi_0,\psi_{-n}). \]
We know that $\calG$ is an affine smooth group scheme over $\calO$
whose generic fiber $\calG_K$ is equal to $\GSp_{2n}$. Furthermore, there
is a left action of $\calG$ on $\bfM^{\rm loc}$.

\subsection{The Kottwitz-Rapoport stratification on $\bfM^{\rm
  loc}_\kappa$}
\label{sec:52}

Let ${\mathcal Fl}$ be the space of chains of $\calO$-lattices in $V$
\[ \pi \calL_0=\calL_{-2n}\subset \dots \subset \calL_{-1}\subset \calL_{0} \]
such that 
\begin{itemize}
\item [(i)] $\calL_{i}/\calL_{i-1}\simeq \kappa$ for each $i$, 
\item [(ii)]  there is a non-degenerate alternating pairing $\psi'$ on
  $\calL_0$ with values in $\calO$ such that $\pi ^m\psi'=\psi$ for some
  $m\in \Z$, and 
\item [(iii)] set $\ol \calL_{-i}:=\calL_{-i}/\calL_{-2n}$, we require
  that the orthogonal complement $\ol \calL_{-i}^\bot$ with respect to
  $\psi'$ is equal to $\ol \calL_{i-2n}$ for all $i$.  
\end{itemize}

Note that a lattice chain $(\calL_\bullet)$ in ${\mathcal Fl}$ is determined
by it members $\calL_{-i}$ for $0\le i\le n$. 

We regard $V$ as the space of column vectors. For $g\in \GSp_{2n}(K)$,
the map 
\[ g \mapsto (\calL_i)=(gL_i) \]
gives a bijection $\GSp_{2n}(K)/I\simeq {\mathcal Fl}$, where $I$ is
the stabilizer of the standard lattice chain. 

Now we restrict to the equi-characteristic case
$\calO=\kappa[[t]]$ and $\pi=t$. 
The space ${\mathcal Fl}$ has a natural ind-scheme structure over
$\kappa$ and is called the {\it affine flag variety} associated to
$\GSp_{2n}$ over $\kappa$. 
For any field extension $\kappa'$ of $\kappa$, we have a
natural bijection
\[ \GSp_{2n}\left(\kappa'((t))\right)/I(\kappa')\simeq {\mathcal
  Fl}(\kappa'), \] 
where $I(\kappa')$ is the stabilizer of the standard lattice chain
$(L_\bullet\otimes \kappa'[[t]])$ base change over $\kappa'[[t]]$. 
 
Let $\scrY$ be the closed subscheme of ${\mathcal Fl}$ consisting of the
lattice chains $\calL_\bullet$ such that 
\[ t L_{-i}\subset \calL_{-i}\subset L_{-i}, \quad 0 \le i \le n, \]
such that $L_0/\calL_0\simeq \kappa^n$. 

The group $I$ acts on the ind-scheme ${\mathcal Fl}$ by the left
translation; 
it leaves the subscheme $\scrY$ invariant. Using the Bruhat-Iwahori
decomposition
\[ \GSp_{2n}(\kappa((t)))=\coprod_{x\in \wt W} IxI, \]
the $I$-orbits are indexed by the extended affine Weyl group $\wt W$
of $\GSp_{2n}$:
\[ {\mathcal Fl}=\coprod_{x\in \wt W} {\mathcal Fl}_x. \]
 
The extended affine Weyl group $\wt W$ is the semi-direct product
$X_*(T)\rtimes W$ of the Weyl group $W$ of $\GSp_{2n}$ and the
cocharacter groups $X_*(T)$, where $T$ is the group of diagonal
matrices in $\GSp_{2n}$. The cocharacter group $X_*(T)$ is 
\[ \{(u_1,\dots, u_{2n})\in \Z^{2n}\, |\,
u_1+u_{2n}=\dots=u_{n}+u_{n+1}, \}. \] 
The Weyl group $W$ is a subgroup of the symmetric group $S_{2n}=W(\GL_{2n})$ 
consisting of elements that commute with the permutation
\[ \theta=(1,2n)(2,2n-1)\dots(n,n+1). \]
We identity the symmetric group $S_{2n}$ with the group of
permutation matrices in $\GL_{2n}$ in the way that

\begin{equation}
  \label{eq:51}
  \forall\, \sigma\in S_{2n}, \text{the corresponding matrix  
  $w_{\sigma}\in \GL_{2n}$ is given by  
  $w_{\sigma}(e_{i})=e_{\sigma(i)}, \forall\,  i$}. 
\end{equation}

We may regard the group $\wt W$ as a subgroup of the group
$\bfA(\R^{2n})$ of affine transformations on the space $\R^{2n}$ of
column vectors. For $\nu\in X_*(T)$, we write $t_\nu$ for the image of
$t$ under $\nu$ in $\GSp_{2n}(\kappa((t)))$, also for the translation
by $\nu$ in 
$\bfA(\R^{2n})$. The element $x=t_\nu w$ then is identified with the
function $x(v)=w\cdot v+\nu$ for $v\in \R^{2n}$. If $x=(x_i)\in
\R^{2n}$, we write $|x|=\sum_i x_i$.   

There is a partial order on the extended affine Weyl group $\wt W$
called the Bruhat order. According to the definition (see
\cite[1.8]{kottwitz-rapoport:alcoves}, for
example), one defines the Bruhat order on the extended affine Weyl
group $\wt W_{\rm der}$ of $\Sp_{2n}$ first, then the partial order on
$\wt W$ is inherited from $\wt W_{\rm der}$ as follows: 
\[ \text{$x\le y$ in $\wt W$ $\iff$
$[x]=[y]$ in $\wt W_{\rm der}\backslash \wt W$ and $1\le y x^{-1}$ in
$\wt W_{\rm der}$.}\] The choice of the Bruhat order on $\wt W_{\rm
der}$ depends on the 
choice of the Borel subgroup of $\Sp_{2n}$. We choose the Borel
subgroup $B$ to be the subgroup of upper triangular matrices in
$\Sp_{2n}$. This also agrees with the choice (following Haines
\cite{haines:clay}) of the lattice chain $L_{-2n}\subset \dots
\subset L_{0}$.    

Let $v_{-2n},\dots, v_0 \in \Z^{2n}$ be the integral vectors
corresponding the lattices $L_{-2n},\dots, L_{0}$. We have 
$v_{-i}=(0^{2n-i},1^{i})$. 
The sequence of integral vectors $v_{-2n},\dots, v_0 \in \Z^{2n}$
forms an alcove in the sense of Kottwitz and Rapoport (see
\cite[Subsection 3.2, 4.2]{kottwitz-rapoport:alcoves}).  Let $\mu=(1^n,0^n)\in X_*(T)$, a dominant
coweight. Following \cite{kottwitz-rapoport:alcoves}, 
we define the sets of 
{\it $\mu$-permissible} and {\it $\mu$-admissible} elements:
\[ {\rm Perm}(\mu):=\{x\in \wt W\,|\, (0,\dots,0) \le
x(v_{-i})-v_{-i}\le (1,\dots,1) \text{ for all $i$ and $|x(0)|=n$}\} \]   
\[ {\rm Adm}(\mu):=\{x\in \wt W\,|\, \text{there is an element $w\in
  W$ such that $x\le t_{w(\mu)}$} \}\]

\begin{prop} \label{51} Notation as above.
\begin{itemize}
\item [(1)] The stratum ${\mathcal Fl}_x$ is contained in $\scrY$ if
  and only if $x\in {\rm Perm}(\mu)$.
\item [(2)] ${\rm Adm}(\mu)={\rm Perm}(\mu)$. 
\end{itemize}
\end{prop}
\begin{proof}
  (1) Let $x\in \wt W\subset \GL_{2n}(\kappa((t)))$. The lattice $xL_{-i}$
  corresponds to the element $x(v_{-i})$ in $\Z^{2n}$. Then the condition
  $tL_{-i}\subset xL_{-i}\subset L_{-i}$ is easily seen to be
  $(0,\dots,0)\le x(v_{-i})-v_{-i}\le (1,\dots,1)$. We also have
  $\dim_\kappa L_0/xL_0=|x(0)|$. Therefore, the statement follows.  

  (2) This is Theorem 4.5 (3) of Kottwitz and Rapoport
      \cite{kottwitz-rapoport:alcoves}. \qed 
\end{proof}

\begin{remark}
The embedding $\sigma \mapsto w_\sigma$ from $S_{2n}$ to
$\GL_{2n}$ in (\ref{eq:51}) does not send the Weyl group $W$ into
$\GSp_{2n}$. In fact for each $\sigma\in W$, there is a unique element
$\epsilon_\sigma={\rm diag}(1,\dots, 1, \epsilon_{\sigma,n+1},\dots,
\epsilon_{\sigma,2n})$ with $\epsilon_{\sigma,i}\in \{\pm 1\}$ such
that $w'_\sigma=w_\sigma \epsilon_\sigma\in \Sp_{2n}$. However, since 
$w_{\sigma}'t_\nu (w_\sigma')^{-1}=w_\sigma t_\nu w_\sigma^{-1}$ and
$t_\nu w_\sigma' L_{-i}= t_\nu w_\sigma L_{-i}$, it won't effect any
results if we choose the presentation of $\wt W$ in $\GL_{2n}$ either 
by $(\nu,\sigma)\mapsto t_\nu w_\sigma$ or by $(\nu,\sigma)\mapsto
t_\nu w_\sigma'$. We make the first choice as it is easier not to deal
with the signs. Another reason for this choice is that the lattice
point in $\Z^{2n}$ corresponding to $t_\nu w_\sigma'L_{-i}=t_\nu
w_\sigma L_{-i}$ is $w_\sigma\cdot v_{-i}+\nu$ not $w_\sigma' \cdot
v_{-i}+\nu$. 
\end{remark} 

To avoid confusing the standard lattice chain that defines the
local model $\bfM^{\rm loc}$ and the lattice chain for the affine flag
variety ${\mathcal Fl}$, we use different notation to distinguish
them. Let 
$\Lambda'_{-i}=\kappa[[t]]^{2n}$ for $0\le i\le 2n$, and
$L'_{i-2n}=<e_1,\dots,e_i,te_{i+1},\dots te_{2n}>$. 
Define $\beta_{-i}'$,
$a_{-i}'$, $\psi_0',\psi'_{-n}$, $\calG'$ as in Subsection
\ref{sec:51}. In particular, $\calG'$ is a smooth affine group scheme
over $\kappa[[t]]$ and one has
\[ \calG'(S)=\Aut(\Lambda'_\bullet\otimes 
\calO_S,\psi_0',\psi_{-n}')
\] 
for any $\kappa[[t]]$-scheme $S$. One also sees that the generic fiber of
$\calG'$ is $\GSp_{2n}$, $\calG'(\kappa[[t]])$ is equal to $I$, and
that the special fiber 
$\calG'_\kappa$ is canonically isomorphic to $\calG_\kappa$.  
 
Using the isomorphism $a_{-i}':\Lambda_{-i}'\simeq L_{-i}'$, we
regard the lattice $\calL_{-i}$, where $(\calL_{-i})$ is a member in
$\scrY$, as a $\kappa[[t]]$-submodule of $\Lambda'_{-i}$
containing $t\Lambda'_{-i}$. 
Then we have an isomorphism $b: \scrY \simeq \bfM^{\rm loc}_\kappa$,
which maps any $\kappa'$-point 
of $\scrY$ to $\bfM^{\rm loc}_\kappa(\kappa')$ by
\[ (\calL_{-i}) \mapsto (\scrF_{-i}), \quad
\scrF_{-i}:=\calL_{-i}/t\Lambda'_{-i} \subset 
\Lambda_{-i}\otimes_{} \kappa'\]
where $\kappa'$ is any field extension of $\kappa$. 

The action of $I$ on the scheme $\scrY$ factors through the quotient
$\calG'(\kappa)$. We also know that the isomorphism $b$ is 
$\calG_\kappa=\calG'_\kappa$-equivariant. Therefore, the
stratification 
\[ \scrY=\coprod_{x\in {\rm Adm}(\mu)} {\mathcal Fl}_x \]
induces a stratification, called the Kottwitz-Rapoport stratification,
on $\bfM^{\rm loc}_\kappa$
\[ \bfM^{\rm loc}_\kappa =\coprod_{x\in {\rm Adm}(\mu)} \bfM^{\rm
  loc}_x \]
so that the $\calG_\kappa$-orbit $\bfM^{\rm loc}_x$ corresponds to
the $I$-orbit ${\mathcal Fl}_x$.

\subsection{Local model diagrams}
\label{sec:53}

Let $\calA'_{g,\Gamma_0(p),N}$ denote the moduli space over
$\Z_{(p)}[\zeta_N]$ that parametrizes equivalence classes of objects 
\[ \begin{CD}
  \ul A_\bullet: \quad \ul A_0 @>\alpha_{1}>> \ul A_1 @>>> 
  \dots @>>>\ul A_{g-1} @>\alpha_{g}>>\ul A_g,
\end{CD}\]
where 
\begin{itemize}
\item each $\ul A_i=(A_i,\lambda_i,\eta_i)$ is a polarized abelian scheme
with a symplectic level-$N$ structure,
\item $\ul A_0$ and $\ul A_g$ are objects in $\calA_{g,1,N}$, and 
\item each $\alpha_i$ is an isogeny of degree $p$ that preserves the level
structures and the polarizations except when $i=1$, and in this case
one has $\alpha_1^*\lambda_1=p\lambda_0$. 
\end{itemize}

There is a natural isomorphism from $\calA_{g,\Gamma_0(p),N}$ to
$\calA_{g,\Gamma_0(p),N}'$ \cite[Propoistion 1.7]{dejong:gamma}. 
We will identify $\calA_{g,\Gamma_0(p),N}$ with
$\calA_{g,\Gamma_0(p),N}'$ via this natural isomorphism. 

Put $n=g$ and $\calO=\Z_p$ in Subsection~\ref{sec:51}. We get a
projective 
$\Z_p$-scheme $\bfM^{\rm loc}$. Let $S$ be a $\Z_p[\zeta_N]$-scheme and
$\ul A_\bullet$ is an object in
$\calA_{g,\Gamma_0(p),N}(S)$. 
A {\it trivialization} $\gamma$ from the de Rham cohomologies $H^1_{\rm
  DR}(A_\bullet/S)$ to $\Lambda_\bullet \otimes \calO_S$ is a
collection of isomorphisms $\gamma_i: H^1_{\rm
  DR}(A_i/S) \to \Lambda_{-i}\otimes \calO_S$ of $\calO_S$-modules  
such that 
\begin{itemize}
\item the diagram 
\[ \begin{CD}
  H^1_{\rm DR} (A_i/S) @>\alpha_{i}^*>> H^1_{\rm DR} (A_{i-1}/S) \\
  @VV\gamma_iV @VV\gamma_{i-1}V \\
  \Lambda_{-i}\otimes \calO_S @>\beta_{-i}>> \Lambda_{-i+1}\otimes
  \calO_S 
\end{CD}\]
commutes for $1\le i\le g$, 
\item if $e_{\lambda_0},e_{\lambda_g}$ are the non-degenerate
  symplectic pairings induced by the principal polarizations
  $\lambda_0, \lambda_g$, then
  $\gamma_i^*\psi_i$ is a scalar multiple of $e_{\lambda_i}$ by
   some element in $\calO_S^\times$ for $i=0,g$. 
\end{itemize}

With the terminology as above, let $\wt \calA_{g,\Gamma_0(p),N}$
denote the moduli space over $\Z_p[\zeta_N]$ that parametrizes
equivalence classes of objects
$(\ul A_\bullet, \gamma)_S$, where
\begin{itemize}
\item $\ul A_\bullet$ is an object in
$\calA_{g,\Gamma_0(p),N}(S)$, and
\item  $\gamma$ is a trivialization from $H^1_{\rm
  DR}(A_\bullet/S)$ to $\Lambda_\bullet \otimes \calO_S$. 
\end{itemize}

The moduli scheme $\wt \calA_{g,\Gamma_0(p),N}$ has two natural
projections $\varphi^{\rm mod}$ and $\varphi^{\rm loc}$.
The morphism 
\[ \varphi^{\rm mod}:\wt \calA_{g,\Gamma_0(p),N}\to
\calA_{g,\Gamma_0(p),N}\otimes \Z_p[\zeta_N]\] 
forgets the
trivialization. The morphism 
\[ \varphi^{\rm
  loc}: \wt \calA_{g,\Gamma_0(p),N}\to \bfM^{\rm loc}\otimes
\Z_p[\zeta_N]\] 
sends an object $(\ul A_\bullet, \gamma)$ to
$(\gamma(\omega_\bullet))$, where $\omega_\bullet=(\omega_i)$ is a
system of $\calO_S$-submodules in the Hodge filtration
\[ 0\to \omega_i\to H^1_{\rm DR}(A_i/S)
\to R^1f_*(\calO_{A_i})\to 0,  \]
and $f:A_i\to S$ is the structure morphism.
Thus, we have the diagram:   
\[ \xymatrix{
 & \wt \calA_{g,\Gamma_0(p),N} \ar[ld]_{\varphi^{\rm mod}}
 \ar[rd]^{\varphi^{\rm loc}} & \\ 
\calA_{g,\Gamma_0(p),N}\otimes \Z_p[\zeta_N] & & 
 \bfM^{\rm loc}\otimes \Z_p[\zeta_N].  
} \]
The moduli scheme $\wt \calA_{g,\Gamma_0(p),N}$ also has a left action
by the group scheme $\calG$. By the work of Rapoport-Zink
 \cite{rapoport-zink}, de Jong \cite{dejong:gamma},
 and Genestier (cf. Remark below Theorem 1.3 of
 \cite{ngo-genestier:alcoves}), we know 
\begin{itemize}
\item [(a)] $\varphi^{\rm mod}$ is a left $\calG$-torsor, and hence it is
  affine and smooth. 
\item [(b)] $\varphi^{\rm loc}$ is smooth, surjective, 
  $\calG$-equivariant, and of relative dimension same as $\varphi^{\rm
    mod}$.  
\end{itemize}
Let 
\[ \calA_{\Gamma_0(p)}:=\calA_{g,\Gamma_0(p),N} \otimes \Fpbar,\quad  
 \wt \calA_{\Gamma_0(p)}:=\wt \calA_{g,\Gamma_0(p),N}\otimes \Fpbar,\quad 
 \bfM^{\rm loc}_{\Fpbar}:= \bfM^{\rm loc}\otimes \Fpbar \]  
be the reduction modulo $p$, respectively. Let $\wt \calA_{\Gamma_0(p),x}$
 be the pre-image of a KR-stratum $\bfM^{\rm loc}_x$. By $(b)$, $\wt
 \calA_{\Gamma_0(p),x}$ is stable under the
 $\calG_{\Fpbar}$-action. Since $\varphi^{\rm mod}$ is a
 $\calG_{\Fpbar}$-torsor, 
the stratification 
\[ \wt \calA _{\Gamma_0(p)}= \coprod_{x\in {\rm Adm(\mu)}}\wt
\calA_{\Gamma_0(p),x} \]
descends to a stratification, called the Kottwitz-Rapoport
stratification, on $\calA_{\Gamma_0(p)}$:
\[  \calA_{\Gamma_0(p)}=\coprod_{x\in {\rm Adm(\mu)}}
\calA_{\Gamma_0(p),x}. \]
Each stratum $\calA_{\Gamma_0(p),x}$ is smooth of dimension same as
$\dim \bfM^{\rm loc}_{x}$, which is the length
$\ell(x)$ of $x$.  

\subsection{g=2} 
\label{sec:54}
We describe the set ${\rm Adm}(\mu)$ of $\mu$-admissible elements and 
the Bruhat order on this set, in the special case where $g=2$. 
The closure $\bar \bfa$ of the base alcove $\bfa$ is the set of
points $u\in \R^{4}$ such that
$u_1+u_{4}=u_2+u_{3}$ and 
\[ 1+u_1\ge u_{4}\ge u_{3}\ge u_{2}. \]
This is obtained from \cite[12.2]{kottwitz-rapoport:alcoves} by applying
the involution $\theta$ since our choice of the standard alcove
$\{v_{-i}\}$ differs from $\{w_i\}$ in
\cite[4.2]{kottwitz-rapoport:alcoves} by the involution $\theta$.  
The simple reflections corresponding to the faces 
\[ u_3=u_4, \quad u_2=u_3, \quad 1+u_1=u_4, \]
are $s_1=(1\,2)(3\,4)$, $s_2=(2\, 3)$ and  
\[ s_0=((-1,0,0,1),(14)): (u_1,u_2,u_3,u_4)\mapsto (u_4-1, u_2,u_3,u_1+1). \]

One checks that $\tau:=((0,0,1,1),(13)(24))\in \wt W$ is the element in ${\rm
  Adm}(\mu)$ that fixes $\bar \bfa$. 
It is not hard to compute the set ${\rm Perm}(\mu)$ from the
  definition. From this and the fact ${\rm Adm}(\mu)={\rm Perm}(\mu)$ 
(Proposition~\ref{51}), we get 
\[
{\rm Adm}(\mu):=
\left\{
\parbox{3in}{
$\tau$, $s_1 \tau$, $s_0 \tau$, $s_2 \tau$,
  $s_0s_1\tau $, $s_0s_2\tau$, $s_1s_2\tau$,   $s_2s_1\tau$,
  $s_1s_0\tau$, $s_0s_1s_0\tau$, $s_1s_0s_2\tau$, $s_2s_1s_2\tau$,
  $s_0s_2s_1\tau $}  
\right\}.  \]
We compute and express these elements $x$ as $(\nu,\sigma)\in
X_*(T)\rtimes W$: 
\begin{eqnarray}
  & \tau = [(0,0,1,1),(13)(24)],
  & s_1 \tau =[(0,0,1,1), (14)(23)], \nonumber \\
  & s_0 \tau=[(0,0,1,1),(1342)], 
  & s_2 \tau=[(0,1,0,1), (1243)], \nonumber\\
  & s_0s_1\tau=[(0,0,1,1),(23)], 
  & s_0s_2\tau=[(0,1,0,1),(12)(34)], \nonumber\\
  & s_1s_2\tau=[(1,0,1,0), (23)], 
  & s_2s_1\tau=[(0,1,0,1), (14)], \label{eq:541} \\
  & s_1s_0\tau=[(0,0,1,1), (14)], 
  & s_0s_1s_0\tau = [(0,0,1,1),(1)], \nonumber\\
  & s_1s_0s_2\tau= [(1,0,1,0),(1)], 
  & s_2s_1s_2\tau= [(1,1,0,0),(1)], \nonumber\\
  & s_0s_2s_1\tau= [(0,1,0,1), (1)]. &  \nonumber
\end{eqnarray}


For a later use, we also express these elements as $t_\nu
w_\sigma$ in $\GL_4(\kappa((t)))$:\\ 

\begin{center}

$ \tau=\begin{pmatrix}
  & & 1 & \\
  & & & 1 \\
  t & & & \\
  & t & & \\
\end{pmatrix}$, \ 
$ s_1\tau=\begin{pmatrix}
  & & & 1 \\
  & & 1 & \\
  & t & & \\
  t & & & \\
\end{pmatrix}$, \\
 
$ s_0 \tau=\begin{pmatrix}
  & 1 & & \\
  & & & 1 \\
  t & & & \\
  & & t & \\
\end{pmatrix}$, \ 
$ s_2\tau =\begin{pmatrix}
  & & 1 & \\
  t & & & \\
  & & & 1 \\
  & t & & \\
\end{pmatrix}$,  \\ 
\end{center}
\begin{eqnarray}
& s_0s_1 \tau=\begin{pmatrix}
  1 & & & \\
  & & 1 & \\
  & t & & \\
  & & & t \\
\end{pmatrix}, \ 
 s_0s_2 \tau=\begin{pmatrix}
  & 1 & & \\
  t & & & \\
  & & & 1 \\
  & & t & \\
\end{pmatrix}, \ 
  s_1s_2 \tau=\begin{pmatrix}
  t & & & \\
  & & 1 & \\
  & t & & \\
  & & & 1 \\
\end{pmatrix},  \nonumber \\
& \label{eq:542}\\
& s_2s_1 \tau=\begin{pmatrix}
  & & & 1 \\
  & t & & \\
  & & 1 & \\
  t & & & \\
\end{pmatrix}, \ 
 s_1s_0 \tau=\begin{pmatrix}
  & & & 1 \\
  & 1 & & \\
  & & t & \\
  t & & & \\
\end{pmatrix}, \
 s_0 s_1s_0\tau = \begin{pmatrix}
  1 & & & \\
  & 1 & & \\
  & & t & \\
  & & & t \\
\end{pmatrix}, \nonumber \\
& \nonumber \\
& s_1s_0s_2\tau=\begin{pmatrix}
  t & & & \\
  & 1 & & \\
  & & t & \\
  & & & 1 \\
\end{pmatrix}, \  
 s_2s_1s_2\tau=\begin{pmatrix}
  t & & & \\
  & t & & \\
  & & 1 & \\
  & & & 1 \\
\end{pmatrix}, \ 
 s_0s_2s_1\tau=\begin{pmatrix}
  1 & & & \\
  & t & & \\
  & & 1 & \\
  & & & t \\
\end{pmatrix}. \nonumber 
\end{eqnarray}

By a result of Ng\^o-Genestier \cite[Theorem
4.1]{ngo-genestier:alcoves}, the points on each 
KR-stratum $\calA_{\Gamma_0(p),x}$ have constant $p$-rank, which we
denote by $p$-rank$(x)$. Moreover, if $x=(\nu,\sigma)$, then
$\text{$p$-rank}(x)$ is given by the formula
\begin{equation}
  \label{eq:54}
 \text{$p$-rank}(x)=\frac{1}{2}\,\# {\rm Fix}(\sigma),\quad {\rm
   where\ \ 
   Fix}(\sigma):=\{\,i; \  \sigma(i)=i\,\}. 
 \end{equation}
Put ${\rm Adm}^i(\mu):=\{x\in {\rm Adm}(\mu);\,
\text{$p$-rank}(x)=i\,\}$. One easily computes the $p$-ranks of
elements in ${\rm Adm}(\mu)$ using (\ref{eq:54}) and gets
\begin{equation}\label{eq:55}
  \begin{split}
  &{\rm Adm}^2(\mu)=\{s_0s_1s_0\tau, s_1s_0s_2\tau, s_2s_1s_2\tau,
  s_0s_2s_1\tau \},\\
 &{\rm Adm}^1(\mu)=\{ s_0s_1\tau , , s_1s_2\tau,   s_2s_1\tau,
  s_1s_0\tau\}, \\
 &{\rm Adm}^0(\mu)=\{\tau, s_1 \tau, s_0 \tau, s_2
  \tau,s_0s_2\tau  \}. \\
\end{split}
\end{equation}
The Bruhat order on ${\rm Adm}(\mu)$ can be described by the following diagram

\begin{figure}[h]
  \begin{center}
\scalebox{0.6}{\includegraphics{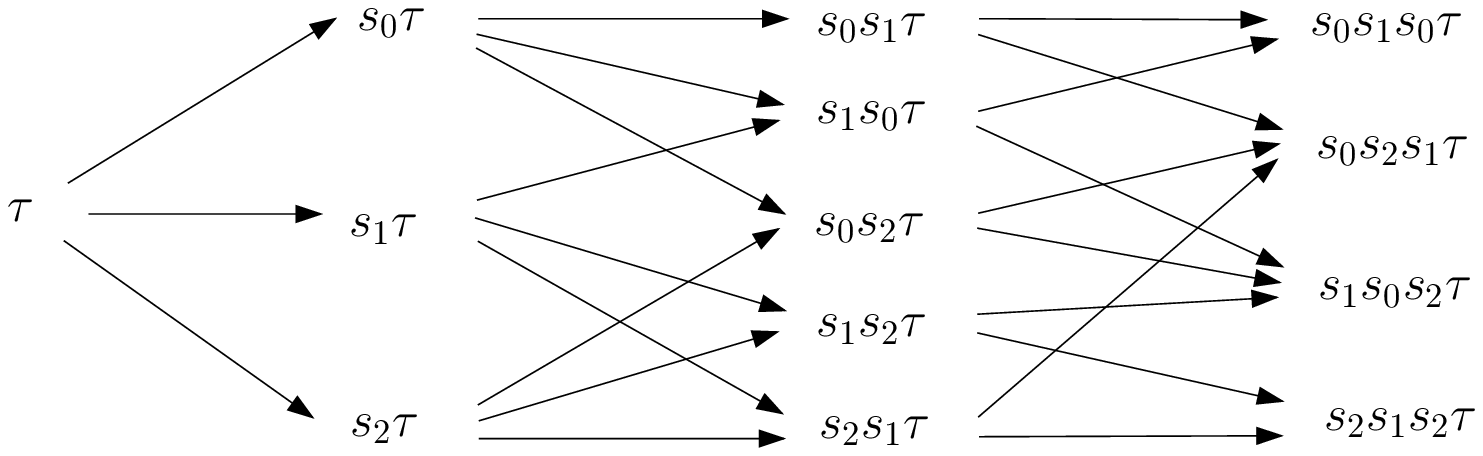}}.  
  \end{center}
\end{figure}
Here $x\to y$ means $\ell(y)=\ell(x)+1$ and $y=sx$ for some reflection $s$
associated an affine root, and $x\le y$ if and only if there exists a
chain 
\[ x=x_1\to x_2\dots \to x_k=y, \]
(cf. \cite[1.1,1.2]{kottwitz-rapoport:alcoves}). This diagram is
obtained  from reading the picture in Haines \cite[Figure 2 in Section
4]{haines:clay} of admissible alcoves for $\GSp_4$

\begin{prop}\label{53} \
{\rm (1)} The supersingular locus
  $\calS_{2,\Gamma_0(p)}$ of
  $\calA_{2,\Gamma_0(p)}$ consists of
  both one-dimensional irreducible components and two-dimensional
  irreducible components.

{\rm (2)} The almost ordinary locus
  $\calA_{2,\Gamma_0(p)}^{1}$ is not dense in the non-ordinary locus
  $\calA_{2,\Gamma_0(p)}^{\rm non-ord}$.

\end{prop}

\begin{proof}
From (\ref{eq:55}) we see that the stratum
$\calA_{\Gamma_0(p),s_0s_2\tau}$ is supersingular and has dimension
two. From the above diagram we read that the stratum
$\calA_{\Gamma_0(p),s_1\tau}$ is not contained in the closure of
$\calA_{\Gamma_0(p),s_0s_2\tau}$. Therefore, we conclude the statement
(1). Since the ordinary locus is dense, the non-ordinary locus has
dimension two. Then the statement (2) follows immediately from (1). 
\end{proof}

\begin{remark}\label{54}
  (1) It follows from Proposition~\ref{53} (2) that for $g\ge 2$ and 
  $1\le f\le g-1$, the $p$-rank-$f$ stratum  $\calA_{\Gamma_0(p)}^f$
  is not dense in the stratum $\calA_{\Gamma_0(p)}^{\le f}\subset
  \calA_{\Gamma_0(p)}$ consisting of points with $p$-rank $\le
  f$. This shows that the collection $\{\calA_{\Gamma_0(p)}^{
  f}\}$ of $p$-rank strata does not form a stratification.  
 
  (2) It follows from Proposition~\ref{53} (1) that for
  $g\ge 2$ and $0\le f\le g-2$, the natural morphism
  $\calA_{\Gamma_0(p)}^f\to \calA^f$ is not finite.   
  
  (3) Tilouine \cite[p.790]{tilouine:coates} examines the
  intersections of 
  four components of the special fiber of the local model. He also
  concludes the same results above. 
\end{remark}

\section{Geometric characterization ($g=2$)}
\label{sec:06}
Let $a=(\ul A_0\stackrel{\alpha}{\to} \ul A_1 \stackrel{\alpha}{\to}
\ul A_2)$ be a point in 
$\calA_{2,\Gamma_0(p)}(k)$, where $k$ is an \ac field of \ch $p$. Then 
$a$ lies in a Kottwitz-Rapoport stratum 
$\calA_{\Gamma_0(p),KR(a)}$ for a unique element $KR(a)$ 
in ${\rm Adm}(\mu)$. We would like to describe
$KR(a)$ from geometric properties of the point $a$. 

Put $M_i:=H^1_{\rm DR}(A_i/k)$ and $\omega_i:=\omega_{A_i}$. We have 
\[ 
\begin{CD}
M_2@>\alpha>> M_1 @>\alpha>> M_0.   
\end{CD}\]
Set $G_0:=\ker (\alpha:A_0\to A_1)$ and $G_1:=\ker (\alpha:A_1\to
A_2)$; they are finite flat group scheme of rank $p$, which is
isomorphic to $\Z/p, \mu_p$, or $\alpha_p$..
From the \dieu theory we know that
\[ \omega_{G_i}=\omega_i/\alpha(\omega_{i+1}),\quad \text{and\ \ }
  \Lie(G_i^t)=M_i/(\omega_i+\alpha(M_{i+1})), \]
where $G_i^t$ is the Cartier dual of $G_i$. Following de Jong
\cite{dejong:gamma}, we define
\[ \sigma_i(a):=\dim \omega_i/\alpha(\omega_{i+1}),\quad \tau_i(a):=\dim
M_i/(\omega_i+\alpha(M_{i+1})).  \]
Clearly we have the following characterization of $G_i$:

\begin{center}
\begin{tabular}{lp{4cm}l}
{\it Table~1.} & 
\end{tabular}\\
\begin{tabular}{|c|c|c|c|}  \hline
$(\sigma_i(a),\tau_i(a))$ & $(0,1)$ & $(1,0)$ & $(1,1)$ \\ \hline
$G_i$ & $\Z/p\Z$ & $\mu_p$ & $\alpha_p$  \\ \hline
\end{tabular} \\ 
\end{center} \ 

When the $p$-rank of $a$ is $\ge 1$, the chain of  the
$p$-divisible groups of $a$ is determined by the invariants 
$(\sigma_i(a), \tau_i(a))$ up to isomorphism.  In particular, the
element $KR(a)$ is determined by the invariants $(\sigma_i(a),
\tau_i(a))$. To describe the correspondence, it suffices to compute
these invariants for the distinguished point $x$ in the stratum
$\bfM^{\rm loc}_x$. 

Recall how to associate a member in $\bfM^{\rm loc}_x$ to an element
  $x=t_\nu w_\sigma$ in ${\rm Adm}(\mu)$ . We first apply $x$ to
  the standard lattice chain and get a lattice chain $(xL'_{-i})_{0\le
  i\le 2}$. It follows from the permissibility that  $tL'_{-i}\subset
  xL'_{-i}\subset L'_{-i}$. 
  . Then there is a lattice $\calL_{-i}$ in
  $\Lambda'_{-i}$ so that its image under the isomorphism
  $\Lambda_{-i}'\simeq L'_{-i}$ is $xL'_{-i}$. This way we associate
  an element $(\calL_{-i}/t\Lambda'_{-i})$ in $\bfM^{\rm loc}_x$ (via
  the isomorphism $b: \scrY\simeq \bfM^{\rm loc}_{\Fpbar}$).
  We use this element to compute the invariants  $(\sigma_i,\tau_i)$.

Below we write
$[L'_0]=(e_1,e_2,e_3,e_4)^t,[L'_{-1}]=(e_1,e_2,e_3,te_4)^t$, and 
$[L'_{-2}]=(e_1,e_2,te_3,te_4)^t$ and write $\ol
\calL_{-i}=\calL_{-i}/t\Lambda'_{-i}$ and $\ol
\Lambda'_{-i}=\Lambda'_{-i}/\calL_{-i}$. Recall that $\beta's$ are the
maps between the lattices $\Lambda'_{-i}$ which correspond to the maps 
$\alpha$ on $M_i$ under a trivialization map $\gamma$. Then the
invariants $(\sigma_i,\tau_i)$ are given by 
\[ \sigma_i=\dim \ol \calL_{-i}/\beta(\ol\calL_{-i-1}), \quad
   \tau_i=\dim \ol \Lambda_{-i}'/\beta(\ol \Lambda_{-i-1}'). \]     

\begin{enumerate}

\item When $x=s_0s_1s_0\tau = \begin{pmatrix}
  1 & & & \\
  & 1 & & \\
  & & t & \\
  & & & t \\
\end{pmatrix}$, we compute 
\[ x[L'_0]=
\begin{pmatrix}
  e_1 \\ e_2 \\te_3 \\ te_4
\end{pmatrix}, \quad x[L'_{-1}]=
\begin{pmatrix}
  e_1 \\ e_2 \\te_3 \\ t^2e_4
\end{pmatrix}, \quad x[L'_{-2}]=
\begin{pmatrix}
  e_1 \\ e_2 \\t^2e_3 \\ t^2e_4
\end{pmatrix}.\]
It follows that 
$ \calL_{-2}=\calL_{-1}=\calL_0=<e_1,e_2,te_3,te_4>$, $\ol
\calL_{-2}=\ol \calL_{-1}=\ol \calL_0=<e_1,e_2>$, and $\ol 
\Lambda'_{-2}=\ol \Lambda'_{-1}=\ol \Lambda'_0=<e_3,e_4>$. 
It follows that   
$\beta(\ol \calL_{-2})=\beta(\ol \calL_{-1})=<e_1,e_2>$, 
$\beta(\ol \Lambda'_{-2})=<e_4>$, and $\beta( \ol
\Lambda'_{-1})=<e_3>$.  This gives
$(\sigma_0,\tau_0)=(0,1)$ and $(\sigma_1,\tau_1)=(0,1)$. 

Note that if $x$ is diagonal, then
$\calL_{-2}=\calL_{-1}=\calL_0$. Therefore, it is enough to 
compute $x[L'_0]$, which is done this way in (2)--(4).

\item When $x=s_0s_2s_1\tau=\begin{pmatrix}
  1 & & & \\
  & t & & \\
  & & 1 & \\
  & & & t \\
\end{pmatrix}$, we compute $x[L'_0]=
\begin{pmatrix}
  e_1 \\ te_2 \\e_3 \\ te_4
\end{pmatrix}$ and obtain
\[ \calL_0=<e_1, te_2,e_3, te_4>,\quad \ol\calL_0=<e_1,e_3>, 
\quad \ol\Lambda'_0=<e_2,e_4>. \]
It follows that  $\beta(\ol \calL_{-1})=<e_1,e_3>$, $\beta(\ol
\calL_{-2})=<e_1>$,  $\beta(\ol \Lambda'_{-1})=<e_2>$, and 
$\beta(\ol \Lambda'_{-2})=<e_2,e_4>$. 
This gives
$(\sigma_0,\tau_0)=(0,1)$ and $(\sigma_1,\tau_1)=(1,0)$.

\item When $x=s_1s_0s_2\tau=\begin{pmatrix}
  t & & & \\
  & 1 & & \\
  & & t & \\
  & & & 1 \\
\end{pmatrix}$, we compute $x[L'_0]=
\begin{pmatrix}
  te_1 \\ e_2 \\ te_3 \\ e_4
\end{pmatrix}$ and obtain
\[ \calL_0=<te_1, e_2,te_3, e_4>,\quad \ol\calL_0=<e_2,e_4>, 
\quad \ol\Lambda'_0=<e_1,e_3>. \]
It follows that  $\beta(\ol \calL_{-1})=<e_2>$, $\beta(\ol
\calL_{-2})=<e_2,e_4>$,   
$\beta(\ol \Lambda'_{-1})=<e_1,e_3>$, and 
$\beta(\ol \Lambda'_{-2})=<e_1>$. 
This gives
$(\sigma_0,\tau_0)=(1,0)$ and $(\sigma_1,\tau_1)=(0,1)$.

\item When $x=s_2s_1s_2\tau=\begin{pmatrix}
  t & & & \\
  & t & & \\
  & & 1 & \\
  & & & 1 \\
\end{pmatrix}$, we compute $x[L'_0]=
\begin{pmatrix}
  te_1 \\ te_2 \\ e_3 \\ e_4
\end{pmatrix}$ and obtain
\[ \calL_0=<te_1, te_2,e_3, e_4>,\quad \ol\calL_0=<e_3,e_4>, \
\ol\Lambda'_0=<e_1,e_2>. \]
It follows that  $\beta(\ol \calL_{-1})=<e_3>$, $\beta(\ol
\calL_{-2})=<e_4>$,  
$\beta(\ol \Lambda'_{-1})=<e_1,e_2>$, and 
$\beta(\ol \Lambda'_{-2})=<e_1,e_2>$. 
This gives
$(\sigma_0,\tau_0)=(1,0)$ and $(\sigma_1,\tau_1)=(1,0)$.

\item When $x=s_0s_1 \tau=\begin{pmatrix}
  1 & & & \\
  & & 1 & \\
  & t & & \\
  & & & t \\
\end{pmatrix}$, we compute $x([L'_0],[L'_{-1}],[L'_{-2}])=
\begin{pmatrix}
  e_1 & e_1 & e_1\\
  te_3 & te_3 & te_3 \\
  e_2 & e_2 & te_2 \\
  te_4& t^2e_4& t^2 e_4
\end{pmatrix}$ and obtain
\[ \calL_0=<e_1, e_2,te_3, te_4>,\quad \ol\calL_0=<e_1,e_2>, \quad
\ol\Lambda'_0=<e_3,e_4>, \]
\[ \calL_{-1}=<e_1, e_2,te_3, te_4>,\quad \ol\calL_{-1}=<e_1,e_2>,
\quad \ol\Lambda'_{-1}=<e_3,e_4>, \]
\[ \calL_{-2}=<e_1, te_2,e_3, te_4>,\quad \ol\calL_{-2}=<e_1,e_3>,
\quad \ol\Lambda'_{-2}=<e_2,e_4>. \]
It follows that  $\beta(\ol \calL_{-1})=<e_1,e_2>$, $\beta(\ol
\calL_{-2})=<e_1>$,  $\beta(\ol \Lambda'_{-1})=<e_3>$, and 
$\beta(\ol \Lambda'_{-2})=<e_4>$. 
This gives
$(\sigma_0,\tau_0)=(0,1)$ and $(\sigma_1,\tau_1)=(1,1)$.

\item When $x=s_1s_2 \tau=\begin{pmatrix}
  t & & & \\
  & & 1 & \\
  & t & & \\
  & & & 1 \\
\end{pmatrix}$, we compute $x([L'_0],[L'_{-1}],[L'_{-2}])=
\begin{pmatrix}
  te_1 & te_1 & te_1\\
  te_3 & te_3 & te_3 \\
  e_2 & e_2 & te_2 \\
  e_4& te_4& te_4
\end{pmatrix}$ and obtain
\[ \calL_0=<te_1, e_2,te_3, e_4>,\quad \ol\calL_0=<e_2,e_4>, \quad
\ol\Lambda'_0=<e_1,e_3>, \]
\[ \calL_{-1}=<te_1, e_2,te_3, e_4>,\quad \ol\calL_{-1}=<e_2,e_4>, \quad
\ol \Lambda'_{-1}=<e_1,e_3>, \]
\[ \calL_{-2}=<te_1, te_2,e_3, e_4>,\quad \ol\calL_{-2}=<e_3,e_4>, \quad
\ol\Lambda'_{-2}=<e_1,e_2>. \]
It follows that  $\beta(\ol \calL_{-1})=<e_2>$, $\beta(\ol \calL_{-2})=<e_4>$, 
$\beta(\ol \Lambda'_{-1})=<e_1,e_3>$, and 
$\beta(\ol \Lambda'_{-2})=<e_1>$. 
This gives
$(\sigma_0,\tau_0)=(1,0)$ and $(\sigma_1,\tau_1)=(1,1)$.

\item When $x=s_2s_1 \tau=\begin{pmatrix}
  & & & 1 \\
  & t & & \\
  & & 1 & \\
  t & & & \\
\end{pmatrix}$, we compute $x([L'_0],[L'_{-1}],[L'_{-2}])=
\begin{pmatrix}
  te_4 & te_4 & te_4\\
  te_2 & te_2 & te_2 \\
  e_3 & e_3 & te_3 \\
  e_1& te_1& te_1
\end{pmatrix}$ and obtain
\[ \calL_0=<e_1, te_2, e_3, te_4>,\quad \ol\calL_0=<e_1,e_3>, \quad
\ol\Lambda'_0=<e_2,e_4>, \]
\[ \calL_{-1}=<te_1, te_2, e_3, e_4>,\quad \ol\calL_{-1}=<e_3,e_4>, \quad 
\ol\Lambda'_{-1}=<e_1,e_2>, \]
\[ \calL_{-2}=<te_1, te_2,e_3, e_4>,\quad \ol\calL_{-2}=<e_3,e_4>, \quad
\ol\Lambda'_{-2}=<e_1,e_2>. \]
It follows that  $\beta(\ol \calL_{-1})=<e_3>$, $\beta(\ol \calL_{-2})=<e_4>$, 
$\beta(\ol \Lambda'_{-1})=<e_2>$, and 
$\beta(\ol \Lambda'_{-2})=<e_1,e_2>$. 
This gives
$(\sigma_0,\tau_0)=(1,1)$ and $(\sigma_1,\tau_1)=(1,0)$. 

\item When $x=s_1s_0 \tau=\begin{pmatrix}
  & & & 1 \\
  & 1 & & \\
  & & t & \\
  t & & & \\
\end{pmatrix}$, we compute $x([L'_0],[L'_{-1}],[L'_{-2}])=
\begin{pmatrix}
  te_4 & te_4 & te_4\\
  e_2 & e_2 & e_2 \\
  te_3 & te_3 & t^2e_3 \\
  e_1& te_1& te_1
\end{pmatrix}$ and obtain
\[ \calL_0=<e_1, e_2,te_3, te_4>,\quad \ol\calL_0=<e_1,e_2>, 
\quad \ol\Lambda'_0=<e_3,e_4>, \]
\[ \calL_{-1}=<te_1, e_2,te_3, e_4>,\quad \ol\calL_{-1}=<e_2,e_4>, \quad
\ol\Lambda'_{-1}=<e_1,e_3>, \]
\[ \calL_{-2}=<te_1, e_2, te_3, e_4>,\quad \ol\calL_{-2}=<e_2,e_4>, \quad
\ol\Lambda'_{-2}=<e_1,e_3>. \]
It follows that  $\beta(\ol \calL_{-1})=<e_2>$, $\beta(\ol \calL_{-2})=<e_2,e_4>$, 
$\beta(\ol \Lambda'_{-1})=<e_3>$, and 
$\beta(\ol \Lambda'_{-2})=<e_1>$. 
This gives
$(\sigma_0,\tau_0)=(1,1)$ and $(\sigma_1,\tau_1)=(0,1)$. \\
\end{enumerate}

We conclude the characterization of $KR(a)$ when $p$-rank of $a$ is
$\ge 1$ in the following table:\\

\begin{center}
\begin{tabular}{lp{12.5cm}l}
{\it Table~2.} & 
\end{tabular}\\
\begin{tabular}{|c|c|c|c|c|c|c|c|c|} \hline
$p$-rank$(a)$ & 2 & 2 & 2 & 2 & 1 & 1& 1 & 1  \\ \hline
$(\sigma_0(a),\tau_0(a))$ & $(0,1)$ & $(0,1)$ & $(1,0)$ & $(1,0)$ 
& $(0,1)$ & $(1,0)$& $(1,1)$ & $(1,1)$ \\ \hline
$(\sigma_1(a),\tau_1(a))$ & $(0,1)$ & $(1,0)$ & $(0,1)$ & $(1,0)$ 
& $(1,1)$ & $(1,1)$& $(1,0)$ & $(0,1)$ \\ \hline

$KR(a)$ & $s_0s_1s_0\tau$  & $s_0s_2s_1\tau$  & $s_1s_0s_2\tau$  
& $s_2s_1s_2\tau$     
& $s_0s_1 \tau$ & $s_1s_2 \tau$ & $s_2s_1 \tau$ & $s_1s_0 \tau$ \\ \hline
\end{tabular} \\ 
\end{center} \ 

When $a$ is supersingular (i.e. any member $A_i$ is a supersingular
abelian variety), the element $KR(a)$ is not determined by the
invariants $(\sigma_i(a),\tau_i(a))$. 
In fact, they are all $(1,1)$, as the schemes $G_i$ 
are isomorphic to $\alpha_p$. We treat this case 
separately.

\section{Geometric characterization ($g=2$): 
the supersingular case}  \label{sec:07}

\subsection{}
\label{sec:71}
We continue with a geometric point $a$ in $\calA_{2,\Gamma_0(p)}(k)$,
and suppose that the point $a$ has $p$-rank $0$. We know that
$(\sigma_i(a),\tau_i(a))=(1,1)$ for $i=0,1$. We define a new invariant
$(\sigma_{02}(a), \tau_{02}(a))$ by
\[ \sigma_{02}(a):=\dim \omega_0/\alpha^2(\omega_2), 
\quad \tau_{02}(a):= \dim
M_0/(\omega_0+\alpha^2(M_2). \]

As in the previous section, we associated a distinguished point in
$\bfM^{\rm loc}_x$ to each
element $x=t_\nu w_\sigma$ in ${\rm Adm}(\mu)$. 
We shall use this point to
calculate the invariant $(\sigma_{02}(a), \tau_{02}(a))$. Below  
$[L'_0], [L'_{-1}], [L'_{-2}], \ol\calL_{-i}$, $\ol \Lambda'_{-i}$ 
are as in the previous section.  

\begin{enumerate}
\item When $x=s_0s_2 \tau=\begin{pmatrix}
  & 1 & &  \\
  t & & & \\
  & & & 1\\
  & & t & \\
\end{pmatrix}$, we compute $x([L'_0],[L'_{-1}],[L'_{-2}])=
\begin{pmatrix}
  te_2 & te_2 & te_2\\
  e_1 & e_1 & e_1 \\
  te_4 & te_4 & t^2e_4 \\
  e_3 & te_3 & te_3
\end{pmatrix}$ and obtain
\[ \calL_0=<e_1, te_2, e_3, te_4>,\quad \ol\calL_0=<e_1,e_3>, \quad
\ol\Lambda'_0=<e_2,e_4>, \]
\[ \calL_{-1}=<e_1, te_2,te_3, e_4>,\quad \ol\calL_{-1}=<e_1,e_4>, 
\quad \ol\Lambda'_{-1}=<e_2,e_3>, \]

\[ \calL_{-2}=<e_1, te_2, e_3, t e_4>,\quad \ol\calL_{-2}=<e_1,e_3>, 
\quad \ol\Lambda'_{-2}=<e_2,e_4>. \]
It follows that  $\beta^2(\ol \calL_{-2})=<e_1>$ and 
$\beta^2 (\ol \Lambda'_{-2})=<e_2>$. 
This gives
$(\sigma_{02},\tau_{02})=(1,1)$. 

\item When $x=s_0 \tau=\begin{pmatrix}
  & 1 & & \\
  & & & 1 \\
  t & & & \\
  & & t & \\
\end{pmatrix}$, we compute $x([L'_0],[L'_{-1}],[L'_{-2}])=
\begin{pmatrix}
  te_3 & te_3 & te_3\\
  e_1 & e_1 & e_1 \\
  te_4 & te_4 & t^2e_4 \\
  e_2 & te_2 & te_2
\end{pmatrix}$ and obtain
\[ \calL_0=<e_1, e_2, te_3, te_4>,\quad \ol\calL_0=<e_1,e_2>, \quad
\ol\Lambda'_0=<e_3,e_4>, \]
\[ \calL_{-1}=<e_1, te_2,te_3, e_4>,\quad \ol\calL_{-1}=<e_1,e_4>, 
\quad \ol\Lambda'_{-1}=<e_2,e_3>, \]
\[ \calL_{-2}=<e_1, te_2, e_3, t e_4>,\quad \ol\calL_{-2}=<e_1,e_3>, 
\quad \ol\Lambda'_{-2}=<e_2,e_4>. \]
It follows that  $\beta^2(\ol \calL_{-2})=<e_1>$ and 
$\beta^2 (\ol \Lambda'_{-2})=0$. 
This gives
$(\sigma_{02},\tau_{02})=(1,2)$. 

\item When $x=s_1 \tau=\begin{pmatrix}
  & 1 & & \\
  & & & 1 \\
  t & & & \\
  & & t & \\
\end{pmatrix}$, we compute $x([L'_0],[L'_{-1}],[L'_{-2}])=
\begin{pmatrix}
  te_4 & te_4 & te_4\\
  te_3 & te_3 & te_3 \\
  e_2 & e_2 & t e_2 \\
  e_1 & te_1 & te_1
\end{pmatrix}$ and obtain
\[ \calL_0=<e_1, e_2, te_3, te_4>,\quad \ol\calL_0=<e_1,e_2>, \quad
\ol\Lambda'_0=<e_3,e_4>, \]
\[ \calL_{-1}=<te_1, e_2,te_3, e_4>,\quad \ol\calL_{-1}=<e_2,e_4>, 
\quad \ol\Lambda'_{-1}=<e_1,e_3>, \]
\[ \calL_{-2}=<te_1, te_2, e_3,  e_4>,\quad \ol\calL_{-2}=<e_3,e_4>, 
\quad \ol\Lambda'_{-2}=<e_1,e_2>. \]
It follows that  $\beta^2(\ol \calL_{-2})=0$ and 
$\beta^2 (\ol \Lambda'_{-2})=0$. 
This gives
$(\sigma_{02},\tau_{02})=(2,2)$.

\item When $x=s_2 \tau=\begin{pmatrix}
  & & 1 & \\
  t & & & \\
  & & & 1 \\
  & t & & \\
\end{pmatrix}$, we compute $x([L'_0],[L'_{-1}],[L'_{-2}])=
\begin{pmatrix}
  te_2 & te_2 & te_2\\
  te_4 & te_4 & te_4 \\
  e_1 & e_1 & t e_1 \\
  e_3 & te_3 & te_3
\end{pmatrix}$ and obtain
\[ \calL_0=<e_1, t e_2, e_3, te_4>,\quad \ol\calL_0=<e_1,e_3>, \quad
\ol\Lambda'_0=<e_2,e_4>, \]
\[ \calL_{-1}=<e_1, te_2,te_3, e_4>,\quad \ol\calL_{-1}=<e_1,e_4>, 
\quad \ol\Lambda'_{-1}=<e_2,e_3>, \]
\[ \calL_{-2}=<te_1, te_2, e_3,  e_4>,\quad \ol\calL_{-2}=<e_3,e_4>, 
\quad \ol\Lambda'_{-2}=<e_1,e_2>. \]
It follows that  $\beta^2(\ol \calL_{-2})=0$ and 
$\beta^2 (\ol \Lambda'_{-2})=<e_2>$. 
This gives
$(\sigma_{02},\tau_{02})=(2,1)$. 

\item When $x=\tau=\begin{pmatrix}
  & & 1 & \\
  & & & 1 \\
  t & & &  \\
  & t & & \\
\end{pmatrix}$, we compute $x([L'_0],[L'_{-1}],[L'_{-2}])=
\begin{pmatrix}
  te_3 & te_3 & te_3\\
  te_4 & te_4 & te_4 \\
  e_1 & e_1 & te_1 \\
  e_2 & te_2 & te_2
\end{pmatrix}$ and obtain
\[ \calL_0=<e_1,  e_2, te_3, te_4>,\quad \ol\calL_0=<e_1,e_2>, \quad
\ol\Lambda'_0=<e_3,e_4>, \]
\[ \calL_{-1}=<e_1, te_2,te_3, e_4>,\quad \ol\calL_{-1}=<e_1,e_4>, 
\quad \ol\Lambda'_{-1}=<e_2,e_3>, \]
\[ \calL_{-2}=<te_1, te_2, e_3,  e_4>,\quad \ol\calL_{-2}=<e_3,e_4>, 
\quad \ol\Lambda'_{-2}=<e_1,e_2>. \]
It follows that  $\beta^2(\ol \calL_{-2})=0$ and 
$\beta^2 (\ol \Lambda'_{-2})=0$. 
This gives
$(\sigma_{02},\tau_{02})=(2,2)$. \\
\end{enumerate}

We conclude the result of our computation for characterizing $KR(a)$
when $p$-rank of $a$ is $0$ in the following table:

\begin{center}
\begin{tabular}{lp{5.8cm}l}
{\it Table~3.} & 
\end{tabular}\\
\begin{tabular}{|c|c|c|c|c|} \hline
$p$-rank$(a)$ & 0 & 0 & 0 & 0 \\ \hline
$(\sigma_0(a),\tau_0(a))$ & $(1,1)$ & $(1,1)$ & $(1,1)$ & $(1,1)$ 
  \\ \hline
$(\sigma_1(a),\tau_1(a))$ & $(1,1)$ & $(1,1)$ & $(1,1)$ & $(1,1)$   
\\ \hline
$(\sigma_{02}(a),\tau_{02}(a))$ & $(1,1)$ & $(1,2)$ & $(2,1)$ & $(2,2)$   
\\ \hline
$KR(a)$ & $s_0s_2\tau$  & $s_0\tau$  & $s_2\tau$  
& $s_1 \tau$, $\tau$  \\ \hline
\end{tabular} \\ 
\end{center} \ 

To distinguish the types $s_1 \tau$ and $\tau$, we need to know a global
description of the supersingular locus.

\subsection{Description of the supersingular locus}
\label{sec:72}
Let $\calA_{2,p,N}$ denote the moduli space of polarized abelian
surfaces of degree $p^2$ with a level-$N$ structure with respect to
$\zeta_N$. Let $\Lambda_{2,1,N}\subset \calA_{2,1,N}\otimes \Fpbar$
denote the subset of superspecial (geometric) points. Let
$\Lambda\subset \calA_{2,p,N}\otimes \Fpbar$ be
the subset of (geometric) points $(A,\lambda,\eta)$ such that $\ker
\lambda\simeq 
\alpha_p\times \alpha_p$. Any member $\ul A$ of $\Lambda$ is
superspecial, as 
$A$ contains $\ker \lambda=\alpha_p\times \alpha_p$. The sets
$\Lambda_{2,1,N}$ and $\Lambda$ are finite, and every member of them
is defined over $\Fpbar$. 

Recall that $\calA_{2,\Gamma_0(p)}$ denotes the reduction modulo $p$
of the Siegel $3$-fold with Iwahori level structure, which
parametrizes equivalence classes of objects $(\ul
A_0\stackrel{\alpha}{\to} \ul A_1\stackrel{\alpha}{\to} \ul A_2)$ in
\ch $p$ with conditions as before (Subsection~\ref{sec:53}). 
Let $\calS_{2,\Gamma_0(p)}\subset
\calA_{2,\Gamma_0(p)}$ denote 
the supersingular locus, the reduced closed subscheme consisting of
supersingular points. Clearly, we have (\ref{eq:55})
\begin{equation}
  \label{eq:721}
  \calS_{2,\Gamma_0(p)}=\coprod_{x\in {\rm Adm}^2(\mu)}
  \calA_{\Gamma_0(p),x} \quad (g=2). 
\end{equation}
For each $\xi=(A_\xi,\lambda_\xi,\eta_\xi) \in \Lambda$, let
$W_\xi\subset \calS_{2,\Gamma_0(p)}$
be the reduced closed subscheme consisting of points $(\ul
A_0 {\to} \ul A_1{\to} \ul A_2)$ such that $\ul A_1\simeq \xi$. 
For each $\gamma=(A_\gamma,\lambda_\gamma,\eta_\gamma)\in
\Lambda_{2,1,N}$, let 
$U_\gamma\subset\calS_{2,\Gamma_0(p)}$ be
the locally closed reduced subscheme
consisting of points  $(\ul A_0 {\to} \ul A_1{\to} \ul A_2)$ such that
$\ul A_0\simeq \gamma$ and $\ul A_1\not\in \Lambda$. Let
$\calS_\gamma$ be the Zariski closure of $U_\gamma$ in 
$\calS_{2,\Gamma_0(p)}$. Clearly, 
$\calS_{\gamma_1}\cap \calS_{\gamma_2}=\emptyset$ 
if $\gamma_1\ne \gamma_2$ and $W_{\xi_1}\cap W_{\xi_2}=\emptyset$ 
if $\xi_1\ne \xi_2$. 

\begin{thm}\label{71} Notation as above. 

{\rm (1)} One has
\[  \calS_{2,\Gamma_0(p)}=\left (\coprod_{\xi\in \Lambda} W_\xi
\right )\cup \left (\coprod_{\gamma\in \Lambda_{2,1,N}} 
S_\gamma\right ) \]
as the union of irreducible components.
Consequently, the supersingular locus has
$|\Lambda|+|\Lambda_{2,1,N}|$ irreducible components. 

{\rm (2)} For each $\xi\in \Lambda$, the subscheme $W_\xi$ is
isomorphic to $\bfP^1\times \bfP^1$ over $\Fpbar$. For each $\gamma\in
\Lambda_{2,1,N}$, the subscheme $S_\gamma$ is isomorphic to $\bfP^1$
over $\Fpbar$. Furthermore, $W_\xi$ and $S_\gamma$ intersects
transversally at most one point. The singular locus 
$\calS_{2,\Gamma_0(p)}^{\rm sing}$ is the intersection
\[  \left (\coprod_{\xi\in \Lambda} W_\xi \right )\cap 
\left (\coprod_{\gamma\in \Lambda_{2,1,N}} S_\gamma\right ). \] 


{\rm (3)} One has $|\calS_{2,\Gamma_0(p)}^{\rm
  sing}|=|\Lambda_{2,1.N}|(p+1)$ and
\begin{equation}
  \label{eq:722}
  \begin{split}
     |\Lambda|& =|\Sp_4(\Z/N\Z)|
\frac{(-1)\zeta(-1)\zeta(-3)}{4} (p^2-1), \\
    |\Lambda_{2,1,N}|&=|\Sp_4(\Z/N\Z)|
 \frac{(-1)\zeta(-1)\zeta(-3)}{4}(p-1)(p^2+1),
 \end{split}
\end{equation}
where $\zeta(s)$ is the Riemann zeta function.
\end{thm}

The proof will be given in Subsection~\ref{sec:76}.

\subsection{}
We use the classical contravariant \dieu theory. 
We refer the reader to
Demazure \cite{demazure} for a basic account of this theory. For a
perfect field $k$ of \ch $p$, write $W:=W(k)$ for the ring of Witt 
vectors over
$k$, and $B(k)$ for the fraction field of $W(k)$. Let $\sigma$ be the
Frobenius
map on $B(k)$. A quasi-polarization on a \dieu module $M$ over $k$
here is a non-degenerate (meaning of non-zero discriminant) 
alternating pairing
\[ \<\, ,\>:M\times M\to B(k), \]
such that $\<Fx,y\>=\<x,Vy\>^\sigma$ for $x, y\in M$ and
$\<M^t,M^t\>\subset W$. Here 
the dual $M^t$ of $M$ is regarded as a \dieu submodule in $M\otimes B(k)$ using
the pairing. A quasi-polarization is called {\it separable} if $M^t=M$.
Any polarized abelian variety $(A,\lambda)$ over $k$ naturally gives rise to a
quasi-polarized \dieu module. The induced quasi-polarization is separable
if and only if $(p,\deg \lambda)=1$. 

Assume that $k$ is an \ac field field of \ch $p$.

\begin{lm}\label{72} \ 

  {\rm (1)} Let $M$ be a separably quasi-polarized superspecial \dieu
  module over $k$ of rank $4$. Then there exists a basis
  $f_1,f_2,f_3,f_4$ for $M$ over $W:=W(k)$ 
  such that
  \[ F f_1=f_3, Ff_3=pf_1, \quad Ff_2=f_4, Ff_4=pf_2 \]
    and the non-zero pairings are
  \[ \<f_1,f_3\>=-\<f_3,f_1\>=\beta_1,\quad 
  \<f_2,f_4\>=-\<f_4,f_2\>=\beta_1,\] 
    where $\beta_1\in W(\F_{p^2})^\times$ with
  $\beta_1^\sigma=-\beta_1$.

  {\rm (2)} Let $\xi$ be a point in $\Lambda$, and let $M_{\xi}$ be
  the \dieu module of $\xi$. Then there is a $W$-basis 
  $e_1, e_2, e_3, e_4$ for $M_\xi$ such that
  \[ Fe_1=e_3,\quad Fe_2=e_4,\quad Fe_3=pe_1, \quad Fe_4=pe_2, \]
  and the non-zero pairings are
  \[  \<e_1,e_2\>=-\<e_2,e_1\>=\frac{1}{p}, \quad
  \<e_3,e_4\>=-\<e_4,e_3\>=1. \] 
\end{lm}

\begin{proof} 
  This is Lemma 4.2 of \cite{yu:ss_siegel}. Statement (1)
  is a special case of Proposition 6.1 of \cite{li-oort}, and the statement
  (2) is deduced from that proposition. \qed
\end{proof}

\begin{lm}\label{73} \

{\rm (1)} Let $(M_0,\<\,,\>_0)$ be a separably quasi-polarized
  supersingular \dieu module of rank $4$ and suppose $a(M_0)=1$. Let
  $M_1:=(F,V)M_0$ and $N$ be the unique \dieu module containing $M_0$
  with $N/M_0=k$. Let $\<\,, \>_1:=\frac{1}{p}\<\,,\>_0$ be 
  the quasi-polarization for $M_1$.  
  Then one has $a(N)=a(M_1)=2$, $VN=M_1$, and 
  $M_1/M_1^t\simeq k\oplus k$ as \dieu modules.

{\rm (2)} Let $(M_1,\<\,,\>_1)$ be a quasi-polarized
  supersingular \dieu module of rank $4$. Suppose that $M_1/M_1^t$ is
  of length $2$, that is, the quasi-polarization has degree $p^2$. 
  \begin{itemize}
  \item [(i)] If $a(M_1)=1$, then letting $M_2:=(F,V)M_1$, one has that 
  $a(M_2)=2$ and $\<\,  ,\>_1$ is a separable quasi-polarization on $M_2$.  
  \item [(ii)] Suppose $(M_1,\<\,,\>_1)$ decomposes as the product of two 
  quasi-polarized \dieu submodules of rank $2$.
  Then there are a unique
  \dieu submodule $M_2$ of $M_1$ with $M_1/M_2=k$ and a unique \dieu
  module $M_0$ containing $M_1$ with $M_0/M_1=k$ so that
  $\<\,,\>_1$ {\rm (}resp.  $p\<\,,\>_1${\rm )} is a separable
  quasi-polarization on 
  $M_2$ {\rm(}resp. $M_0${\rm)}. 
  \item[(iii)] Suppose $M_1/M_1^t\simeq k\oplus k$ as \dieu modules. Let
  $M_2\subset M_1$ be any \dieu submodule with $M_1/M_2=k$, and
  $M_0\supset M_1$ be any \dieu overmodule with $M_0/M_1=k$. Then
  $\<\,,\>_1$ {\rm (}resp.  $p\<\,,\>_1${\rm )} is a separable 
  quasi-polarization on $M_2$ {\rm(}resp. $M_0${\rm)}.
  \end{itemize}
\end{lm}
This is well-known; the proof is elementary and omitted. We remark
that \dieu lattices in a supersingular four-dimensional polarized
isocrystal are also classified in C. Kaiser \cite[Section
3]{kaiser:thesis}. 

\subsection{} \label{sec:74} Let $(A_0,\lambda_0)$ be a superspecial
principally polarized abelian surface and $(M_0,\<\,,\>_0)$ be the
associated 
\dieu module. Let $\varphi:(A_0,\lambda_0)\to (A,\lambda)$ be an
isogeny of degree $p$ with $\varphi^* \lambda=p\,\lambda_0$.
Write $(M,\<\, ,\>)$ for the \dieu module of $(A,\lambda)$. Choose a basis
$f_1,f_2,f_3,f_4$ for $M_0$ as in Lemma~\ref{72}. Put $M_2:=(F,V)M_0=VM_0$. 
We have the inclusions
  \[ M_2\subset M \subset M_0. \]
Modulo $M_2$, a module $M$ corresponds a one-dimensional subspace
$M/M_2$ in $M_0/M_2$. As $M_0/M_2=k<f_1,f_2>$, the subspace $M/M_2$
has the form
\[ k<af_1+bf_2>, \quad [a:b]\in \bfP^1(k). \]
Let $\ol M_0:=M_0/pM_0$, and let
\[ \<\, ,\>_0:\ol M_0\times \ol M_0\to k. \]
be the induced perfect pairing.

\begin{lm}\label{74}
Notation as above, the following conditions are equivalent 
\begin{itemize}
\item [(a)] $\ker \lambda\simeq \alpha_p\times \alpha_p$.
\item [(b)] $\<\ol M, F\ol M\>_0=0$, where $\ol M:=M/pM_0$.
\item [(c)] $\<\ol M, V\ol M\>_0=0$.
\item [(d)] The corresponding point $[a:b]$ satisfies $a^{p+1}+b^{p+1}=0$  
\end{itemize}
\end{lm}
\begin{proof}
One has $\ol M=k<f_1',f_3,f_4>$ with $f'_1=af_1+bf_2$. It is easy to
see that 
\[ \<\ol M, F\ol M\>_0=0\iff \<f'_1, Ff'_1\>_0=0\iff  a^{p+1}+b^{p+1}=0, \]
and 
\[ \<\ol M, V\ol M\>_0=0\iff \<f'_1, Vf'_1\>_0=0\iff  a^{p+1}+b^{p+1}=0. \]
This shows that the conditions (b), (c) and (d) are equivalent.

Since $\varphi^*\lambda=p\lambda_0$, we have 
$\<\, ,\>=\frac{1}{p}\<\, ,\>_0$. The \dieu module $M(\ker \lambda)$
of the subgroup $\ker \lambda $ is equal to $M/M^t$. Hence
the condition (a) is equivalent to that $F$ and $V$ vanish on $M(\ker
\lambda)=M/M^t$. On the other hand, the subspace $\ol {M^t}:=M^t/pM_0$
is equal to $\ol M^{\bot}$ with respect to $\<\, ,\>_0$. It follows 
that (a) is equivalent to the conditions (b) and (c).   \qed
\end{proof}

It follows from Lemma~\ref{74} that there are $p+1$ isogenies
$\varphi$ so that  $\ker \lambda\simeq \alpha_p\times \alpha_p$. 
Conversely, fix a polarized superspecial abelian surface
$(A,\lambda)$ such that $\ker \lambda \simeq \alpha_p\times \alpha_p$.
Then there are $p^2+1$ degree-$p$ isogenies $\varphi:(A_0,\lambda_0)\to
(A,\lambda)$ such that $A_0$ is superspecial and
$\varphi^*\lambda=p\,\lambda_0$. Indeed, each isogeny $\varphi$ always
has the property $\varphi^*\lambda=p\,\lambda_0$ for a principal
polarization $\lambda_0$ (Lemma~\ref{73} (iii)), and there are
$|\bfP^1(\F_{p^2})|$ isogenies with $A_0$ superspecial.

\subsection{}
\label{sec:75}
Let $\calA_P$ be the moduli space of isogenies  
$\alpha:\ul A_0\to \ul A_1$ of degree $p$, where $\ul A_0$ 
is an object in $\calA_{2,1,N}$ and $\ul A_1$ is an object 
in $\calA_{2,p,N}$ such that
$\alpha^*\lambda_1=p \lambda_0$ and $\alpha_* \eta_0=\eta_1$. 
Let $\calS_P\subset \calA_P\otimes \Fpbar$ be the supersingular locus,
the reduced closed subscheme consisting of supersingular points. 
For each $\xi=(A_\xi, \lambda_\xi,\eta_\xi)\in \Lambda$, let
$V_\xi\subset \calS_P$ be the closed 
subvariety consisting of the isogenies 
$\alpha:\ul  A_0\to \ul A_1$ such that $\ul A_1 = \xi$. 
For each $\gamma=(A_\gamma,\lambda_\gamma,\eta_\gamma) 
\in \Lambda_{2,1,N}$, let $S'_\gamma\subset \calS_P$ be the
closed subvariety consisting of the isogenies 
$\alpha:\ul  A_0\to \ul A_1$ such that $\ul A_0 = \gamma$.

It is known that the varieties $V_\xi$ and $S'_\gamma$ are isomorphic
to $\bfP^1$ over $\Fpbar$ (cf. \cite{katsura-oort:surface}). 
We also know 
(\cite[Proposition 4.5]{yu:ss_siegel}) that 
\[ \calS_P=\left (\coprod_{\xi\in \Lambda} V_\xi \right )
\cup \left (\coprod_{\gamma\in \Lambda_{2,1,N}} 
S'_\gamma\right
) \]  
as the union of irreducible components.

\def\pr{{\rm pr}}

Let $\pr:\calS_{2,\Gamma_0(p)}\to \calS_P$ be the natural projection.

\subsection{Proof of Theorem~\ref{71}}
\label{sec:76}
(1) It is easy to see that
\[ \calS_{2,\Gamma_0(p)}=\left (\coprod_{\xi\in \Lambda} W_\xi \right )
\coprod \left (\coprod_{\gamma\in \Lambda_{2,1,N}} 
U_\gamma\right ). \]
The statement follows from this.

(2) Clearly we have
\[ W_\xi \simeq V_\xi\times V'_\xi\simeq \bfP^1\times \bfP^1\quad
\text{(over $\Fpbar$)}, \]
where $V'_\xi$ is the variety parameterizing isogenies
$\alpha:{\xi}\to \ul A_2$ of degree $p$ with $\ul
A_2$ in $\calA_{2,1,N}\otimes \Fpbar$ satisfying $\alpha^*
\lambda_2=\lambda_\xi$ and $\alpha_*\eta_\xi=\eta_2$. This completes
the first assertion.  

Let $a=(\ul A_0\stackrel{\alpha}{\to} \ul A_1 \stackrel{\alpha}{\to}
\ul A_2)$ be a point in 
$U_\gamma(k)$. Since $\ul A_2$ is determined by $\ul A_1$
(cf. Lemma~\ref{73} (i) (ii)), the projection $\pr$ induces an
isomorphism
\[ \pr: U_\gamma\stackrel{\sim}{\to} \pr(U_\gamma)\subset S'_\gamma. \]
As $U_\gamma$ is dense in $S_\gamma$ and $S'_\gamma$ is proper,
$\pr(S_\gamma)\subset S'_\gamma$.
Since $\pr$ is proper and $S'_\gamma$ is a smooth curve, 
the section $s:\pr(U_\gamma)\to U_\gamma$ extends uniquely to a
section $s: S'_\gamma\to S_\gamma$. This shows $\pr: S_\gamma\simeq
S'_\gamma$, and hence $S_\gamma\simeq \bfP^1$ over $\Fpbar$.

A component $S_\gamma$ meets a component $W_\xi$ if and only if
$S'_\gamma$ and $\pr(W_\xi)$ meet. Since $S'_\gamma$ and $V_\xi$
meet transversally at most one point \cite[Proposition
4.5]{yu:ss_siegel}, the components $S_\gamma$ and $W_\xi$ meet
transversally at most one point. Since any irreducible component of
$\calS_{2,\Gamma_0(p)}$ is smooth, the singularity occurs only 
at the intersection of components $S_\gamma$ and $W_\xi$.    
 
(3) We know that $S'_\gamma$ contains $p+1$ points with $\ul A_1\in
    \Lambda$ (Subsection~\ref{sec:74}). 
    Each component $S_\gamma$ meets $p+1$ components
    of the form $W_\xi$, and hence has $p+1$ singular points. This
    proves the first part.

  The result (\ref{eq:722}) is due to Katsura and Oort
  \cite[Theorem 5.1, Theorem 5.3]{katsura-oort:surface} in a slightly
  different form. For another proof (using a mass formula due to
  Ekedahl \cite{ekedahl:ss} and some others), see \cite[Corollary
  3.3, Corollary 4.6]{yu:ss_siegel}.   \qed 

\subsection{}
\label{sec:77}
It follows from the description of the supersingular locus that
\begin{itemize}
\item [(i)] The closure of the stratum
  $\calA_{\Gamma_0(p),s_0s_2\tau}$ is $\coprod_{\xi\in
  \Lambda} W_\xi$. 
\item [(ii)] The stratum $\calA_{\Gamma_0(p),s_1\tau}$ is 
$\coprod_{\gamma\in \Lambda_{2,1,N}}  U_\gamma$, as this is the
complement of the closure $\ol{\calA}_{\Gamma_0(p),s_0s_2\tau}$ in
the supersingular locus $\calS_{2,\Gamma_0(p)}$.
\item [(iii)] The minimal stratum $\calA_{\Gamma_0(p), \tau}$ is
\[  \left (\coprod_{\xi\in \Lambda} W_\xi \right )\cap 
\left (\coprod_{\gamma\in \Lambda_{2,1,N}} S_\gamma\right )
=\calS_{2,\Gamma_0(p)}^{\rm sing}. \] 
\end{itemize}

We compute the locus $\ol {\calA}_{\Gamma_0(p),s_0\tau}$, the closure
of the stratum $\calA_{\Gamma_0(p),s_0\tau}$. This is the disjoint
union of the subvarieties in the component $W_\xi$, for $\xi \in
\Lambda$, defined by the closed condition $\tau_{02}=2$ (Table 3.). 

Let $\ul A_1=\xi\in \Lambda$ and $(M_1,\<\,,\>_1)$ be the \dieu module of
$\ul A_1$. Choose a basis $e_1,e_2,e_3,e_4$ for $M_1$ as in
Lemma~\ref{72}. Let $M_2\subset M_1\subset M_0$ be a chain of \dieu
modules with $M_0/M_1\simeq M_1/M_2\simeq k$. As
$M_1/VM_1=k<e_1,e_2>$, the subspace $M_2/VM_1$ has the form
\[ k<ae_1+b e_2>, \quad [a:b]\in \bfP^1(k). \]
As $V^{-1}M_1/M_1=k<\frac{1}{p}e_3, \frac{1}{p}e_4>$, the subspace
$M_0/M_1$ has the form 
\[ k<c\frac{1}{p}e_3+d \frac{1}{p}e_4>, \quad [c:d]\in \bfP^1(k). \]
Use this as coordinates for $W_{\xi}$, we get and fix an isomorphism 
$\Phi:W_\xi\simeq \bfP^1\times \bfP^1$. Let $A, B, C, D$ be lifts in
$W$ of $a,b,c,d$, respectively. We have
\[ M_2=<Ae_1+B e_2, pe_1, pe_2,e_3, e_4>,\ \text{ and }\ M_0=<e_1,
e_2,e_3, e_4,C \frac{1}{p}e_3+D \frac{1}{p}e_4 >. \] 

The condition $\tau_{02}=2$ says that in $\ol M_0:=M_0/pM_0$, the
subspaces $\ol{VM_0}$ and $\ol M_2$ generates a two-dimensional
subspace. As both have dimension two, the condition means $VM_0=M_2$. 
One has
\[ VM_0=<C^{\sigma^{-1}}e_1+D^{\sigma^{-1}} e_2, pe_1, pe_2,e_3,
e_4>. \]
As both submodules contain $VM_1$, modulo $VM_1$, we get
\[ <c^{p^{-1}} e_1+ d^{p^{-1}} e_4>=<ae_1+be_2>. \]
This gives the equation $a^pd-b^p c=0$. We have shown that
\[ \Phi(\ol {\calA}_{\Gamma_0(p),s_0\tau}\cap W_\xi)=
\{\,([a:b],[a^p:b^p])\in \bfP^1\times \bfP^1\,;\, [a:b]\in
\bfP^1\,\}\simeq \bfP^1.  \] 
This is the graph of the relative Frobenius morphism
$F_{\bfP^1/{\F_p}}:\bfP^1\to \bfP^1$. 
We carry out the similar computation and get
\[ \Phi(\ol {\calA}_{\Gamma_0(p),s_2\tau}\cap W_\xi)=
\{\,([c^p:d^p],[c:d])\in \bfP^1\times \bfP^1\,;\, [c:d]\in
\bfP^1\,\}\simeq \bfP^1.  \]  
This is the transpose of the graph of the relative Frobenius morphism.

We summarize the results as follows.
\begin{prop}\label{75}
  We have 
\[ \ol{\calA}_{\Gamma_0(p),s_0s_2\tau}=\coprod_{\xi\in \Lambda} W_\xi,
\quad \ol{\calA}_{\Gamma_0(p),s_1\tau}=\coprod_{\gamma\in
  \Lambda_{2,1,N}}  S_\gamma, \]
\[  \ol{\calA}_{\Gamma_0(p),s_0\tau}\simeq \coprod_{\xi\in \Lambda}
  \bfP^1,\quad \ol{\calA}_{\Gamma_0(p),s_2\tau}\simeq 
  \coprod_{\xi\in \Lambda}\bfP^1, \quad
  \calA_{\Gamma_0(p),\tau}=\calS_{2,\Gamma_0(p)}^{\rm sing}. \]
Consequently, we have
\begin{itemize}
\item [(i)] The stratum $\calA_{\Gamma_0(p),s_0s_2\tau}$ has $|\Lambda|$
  irreducible components.
\item [(ii)] The stratum $\calA_{\Gamma_0(p),s_1\tau}$ has
  $|\Lambda_{2,1,N}|$ 
  irreducible components.
\item [(iii)] The stratum $\calA_{\Gamma_0(p),s_0\tau}$ has $|\Lambda|$
  irreducible components.
\item [(iv)] The stratum $\calA_{\Gamma_0(p),s_2\tau}$ has $|\Lambda|$
  irreducible components.
\item [(v)] The stratum $\calA_{\Gamma_0(p),\tau}$ consists of
  $|\Lambda_{2,1,N} | (p+1)$ points. 
\end{itemize}
\end{prop}

Proposition~\ref{75}, Theorem~\ref{41}, and Proposition 2.1 of
\cite{yu:gamma} answer the question on irreducible components of each
Kottwitz-Rapoport stratum in the moduli space 
$\calA_{2,\Gamma_0(p)}$. \\

We end this paper with the following criterion to distinguish the
types $s_1\tau$ and $\tau$. This finishes our geometric
characterization of Kottwitz-Rapoport strata for $g=2$.  
\begin{lm}
  Let $a=(\ul A_0\stackrel{\alpha}{\to} \ul A_1
\stackrel{\alpha}{\to} \ul A_2)$ be a point in 
$\ol{\calA}_{\Gamma_0(p),s_1\tau}(k)$, and let 
$\ol M_2\stackrel{\alpha}{\to} \ol M_1 \stackrel{\alpha}{\to} 
\ol M_0$ be the chain of the associated de
Rham cohomologies. Let $\omega_i:=\omega_{A_i}\subset \ol M_i$ be the
Hodge subspace. Then the point $a$ lies in $ \calA_{\Gamma_0(p),\tau}$
if and only if the condition $\<\alpha(\ol
M_1),\alpha(\omega_1)\>_0=0$ holds.  
\end{lm}
\begin{proof}
  It follows from Theorem~\ref{71} that $a$ lies in  $
  \calA_{\Gamma_0(p),\tau}$ if and only if the object $\ul A_1$ lies
  in $\Lambda$ (Subsection~\ref{sec:72}). 
  The statement then follows from Lemma~\ref{74}, as one has
  $\omega_1=V\ol M_1$. \qed 
\end{proof}

\begin{thank}
  The author would like to express his appreciation to C.-L.~ Chai,
  M.-T. Chuan and J.~Tilouine for helpful discussions. Part of the
  manuscript is prepared during the author's 
  stay at MPIM in Bonn. 
  He wishes to thank the Institute for kind hospitality and 
  excellent working environment. Finally, he thanks the referee for
  careful reading and helpful comments. 
 
\end{thank}

\end{document}